\title{Uniformly cross intersecting families}
\author{{Noga Alon\thanks{School of Mathematics, Institute for Advanced Study, Princeton, NJ 08540, USA, and Raymond and Beverly Sackler
Faculty of Exact Sciences, Tel Aviv University, Tel Aviv, 69978,
Israel. Email: nogaa@tau.ac.il. Research supported in part by the
Israel Science Foundation, by a USA-Israeli BSF grant, by NSF grant
CCR-0324906, by a Wolfensohn fund and by the State of New Jersey.}}
\quad {Eyal Lubetzky
\thanks{ School of Computer Science, Raymond and Beverly
Sackler Faculty of Exact Sciences, Tel Aviv University, Tel Aviv,
69978, Israel. Email: lubetzky@tau.ac.il. Research partially
supported by a Charles Clore Foundation Fellowship.}}}

\documentclass [11pt]{article}
\usepackage[]{amsmath,amssymb,amsfonts,latexsym,amsthm,enumerate}
\usepackage[]{graphicx}

\newtheorem{theorem}{Theorem}[section]
\newtheorem{lemma}[theorem]{Lemma}
\newtheorem{claim}[theorem]{Claim}
\newtheorem*{definition}{Definition}

\newtheorem{proposition}[theorem]{Proposition}

\renewcommand{\epsilon}{\varepsilon}

\newcommand{\calA}{\mathcal{A}}
\newcommand{\calB}{\mathcal{B}}

\newtheoremstyle{upright}%
        {8pt plus2pt minus4pt}%
        {8pt plus2pt minus4pt}%
        {\upshape}%
        {}%
        {\bfseries}%
        {:}%
        {1em}%
        {}%
\theoremstyle{upright}
\newtheorem{remark}[theorem]{Remark}

\newcommand{\ignore}[1]{}

\begin{document}
\maketitle

\begin{abstract}
Let $\calA$ and $\calB$ denote two families of subsets of an
$n$-element set. The pair $(\mathcal{A},\mathcal{B})$ is said to be
$\ell$-cross-intersecting iff $|A\cap B| = \ell$ for all
$A\in\mathcal{A}$ and $B\in\mathcal{B}$. Denote by $P_\ell(n)$ the
maximum value of $|\mathcal{A}||\mathcal{B}|$ over all such pairs.
The best known upper bound on $P_\ell(n)$ is $\Theta(2^n)$, by
Frankl and R\"{o}dl. For a lower bound, Ahlswede, Cai and Zhang
showed, for all $n \geq 2\ell$, a simple construction of an
$\ell$-cross-intersecting pair $(\mathcal{A},\mathcal{B})$ with
$|\mathcal{A}||\mathcal{B}| =
\binom{2\ell}{\ell}2^{n-2\ell}=\Theta(2^n/\sqrt{\ell})$, and
conjectured that this is best possible. Consequently, Sgall asked
whether or not $P_\ell(n)$ decreases with $\ell$.

In this paper, we confirm the above conjecture of Ahlswede et al.
for any sufficiently large $\ell$, implying a positive answer to the
above question of Sgall as well. By analyzing the linear spaces of
the characteristic vectors of $\mathcal{A},\mathcal{B}$ over
$\mathbb{R}$, we show that there exists some $\ell_0>0$, such that
$P_\ell(n) \leq \binom{2\ell}{\ell}2^{n-2\ell}$ for all $\ell \geq
\ell_0$. Furthermore, we determine the precise structure of all the
pairs of families which attain this maximum.
\end{abstract}

\section{Introduction}\label{sec::intro}
Let $\calA$ and $\calB$ denote two families of subsets of an
$n$-element set. We say that the pair $(\calA,\calB)$ is
$\ell$-cross-intersecting iff $|A\cap B| = \ell$ for all $
A\in\calA$ and $B\in\calB$. Let $P_\ell(n)$ denote the maximum
possible value of $|\calA||\calB|$ over all
$\ell$-cross-intersecting pairs $(\calA,\calB)$. We are interested
in finding the precise value of $P_\ell(n)$, and in characterizing
all the extremal pairs $\calA,\calB$ which achieve this maximum.

The study of the maximal size of a single family of sets
$\mathcal{F} \subset 2^{[n]}$, with specified pairwise intersections
of its members, has received a considerable amount of attention over
the years. For instance, the Erd\H{o}s-Ko-Rado Theorem
\cite{ErdosKoRado}, one of the most fundamental theorems in
Combinatorial Set Theory, gives a tight upper bound $|\mathcal{F}|
\leq \binom{n-t}{k-t}$ in case $|F\cap F'| \geq t$ for all
$F,F'\in\mathcal{F}$, $|F|=k$ for all $F\in\mathcal{F}$ and $n$ is
sufficiently large. The case where there is no restriction on the
size of the sets of $\mathcal{F}$ is treated by Katona's Theorem
\cite{Katona}. In both cases, there is a unique (up to a relabeling
of the elements of $[n]$) family of sets which achieves the upper
bound. For further results of this nature, see, e.g, \cite{Frankl},
\cite{FranklFuredi}, \cite{FranklWilson}, \cite{RCW}, as well as
\cite{BabaiFrankl}.

A well known conjecture of Erd\H{o}s \cite{Erdos2} stated that if
$\mathcal{F}\subset 2^{[n]}$ is a family satisfying $|F\cap F'| \neq
\lfloor \frac{n}{4} \rfloor$ for all $F,F'\in\mathcal{F}$, then
$|\mathcal{F}|<(2-\epsilon)^n$ for some $\epsilon > 0$. This was
proved by Frankl and R\"{o}dl \cite{FranklRodl}, by considering the
corresponding variant on two families: it is shown in
\cite{FranklRodl}, that if $\calA,\calB\subset 2^{[n]}$ and $|A\cap
B| \neq l$, where $\eta n \leq l \leq (\frac{1}{2}-\eta)n$ for some
$\eta < \frac{1}{4}$, then $|\calA||\calB| \leq
(4-\epsilon(\eta))^n$. The authors of \cite{FranklRodl} studied
several additional problems related to cross-intersections of two
families of sets, and among their results, they provided the
following upper bound on $P_\ell(n)$, which was later reproved in
\cite{ACZ}:
\begin{equation}\label{eq-frankl-rodl} \left\{\begin{array}{ll}
P_0(n) \leq 2^n & \ \\
P_\ell(n) \leq 2^{n-1} & \mbox{for }\ell \geq 1
\end{array}\right.~.\end{equation}
The argument which gives the upper bound of $2^n$ is simple:
consider the characteristic vectors of the sets in $\calA,\calB$ as
vectors in $\mathbb{Z}_2^n$. Notice that the intersection of two
sets is equal to the inner product of the two corresponding vectors
modulo $2$. Therefore, if $\ell$ is even, then the families
$\calA,\calB$ belong to two orthogonal linear spaces, giving
$|\calA||\calB| \leq 2^n$. Otherwise, we may add an additional
coordinate of $1$ to all vectors, and repeat (carefully) the above
argument, gaining a slight improvement: $|\calA||\calB| \leq
2^{n-1}$. Similar ideas are used to show that the upper bound
$2^{n-1}$ holds for even values of $\ell > 0$ as well, by performing
the analysis over $GF(p)$ for some prime $p > 2$ instead of over
$\mathbb{Z}_2$.

As part of their study of questions in Coding Theory, Ahlswede, Cai
and Zhang \cite{ACZ} gave the following simple construction of an
$\ell$-cross-intersecting pair: for $n \geq 2\ell$, let $\calA$
contain a single $2\ell$-element set, $A$, and let $\calB$ contain
all the sets which contain precisely $\ell$ elements of $A$. This
gives:
\begin{equation}\label{eq-P-ell-lower-bound}
|\calA||\calB| = \binom{2\ell}{\ell}2^{n-2\ell} =
(1+o(1))\frac{2^n}{ \sqrt{\pi\ell}}~, \end{equation} where the
$o(1)$-term tends to $0$ as $\ell\to\infty$. The upper bound
\eqref{eq-frankl-rodl} implies that this construction achieves the
maximum of $P_\ell(n)$ for $\ell\in\{0,1\}$, and the authors of
\cite{ACZ} conjectured that this in fact holds for all $\ell$.

As the upper bound \eqref{eq-frankl-rodl} is independent of $\ell$,
compared to the above lower bound of $\Theta(2^n/\sqrt{\ell}$),
Sgall \cite{Sgall} asked whether or not $P_\ell(n)$ is bounded from
above by some decreasing function of $\ell$. One of the motivations
of \cite{Sgall} was a relation between problems of restricted
cross-intersections of two families of sets and problems in
Communication Complexity; see \cite{Sgall} for more details.

In \cite{KeevashSudakov}, the authors verified the above conjecture
of \cite{ACZ} for the case $\ell=2$, by showing that $P_2(n) \leq 3
\cdot 2^{n-3}$. However, for any $\ell > 2$ the best known upper
bound on $P_\ell(n)$ remained $2^{n-1}$.

The following theorem confirms the above conjecture of \cite{ACZ}
for all sufficiently large values of $\ell$, and thus provides also
a positive answer to the above question of Sgall.
\begin{theorem}\label{thm-1}
There exists some $\ell_0 > 0$ such that, for all $\ell \geq
\ell_0$, every $\ell$-cross-intersecting pair $\calA,\calB \subset
2^{[n]}$ satisfies:
\begin{equation}\label{eq-final-cross-bound}
|\calA| |\calB| \leq \binom{2\ell}{\ell}2^{n-2\ell}~.
\end{equation}
Furthermore, if $|\calA||\calB|=\binom{2\ell}{\ell}2^{n-\ell}$, then
there exists some choice of parameters $\kappa,\tau,n'$:
\begin{equation}\label{eq-opt-pair-params}
\begin{array}
  {l}
  \kappa\in\{2\ell-1,2\ell\}~,~\tau \in \{0,\ldots,\kappa\}~,\\
  \kappa+\tau \leq n' \leq n,
\end{array}\end{equation}
such that, up to a relabeling of the elements of $[n]$ and swapping
$\calA,\calB$, the following holds: \begin{equation}
\label{eq-opt-pair}
\begin{array}{lll}
        \mathcal{A} &= \Bigg\{
    \displaystyle{\bigcup_{T \in J} T} ~:~J
    \subset\bigg\{\begin{array}{c}
           \{1,\kappa+1\},\ldots,\{\tau,\kappa+\tau\},\\
           \{\tau+1\},\ldots,\{\kappa\}
    \end{array}\bigg\}~,~
    |J|=\ell
    &\Bigg\}\times 2^X~,\\
    \noalign{\medskip}
        \mathcal{B}&=\Bigg\{L \cup \{\tau+1,\ldots,\kappa\} :
        \begin{array}{l}
        L \subset \{1,\ldots,\tau,\kappa+1,\ldots,\kappa+\tau\}\\
        |L\cap \{i,\kappa+i\}| = 1 \mbox{ for all
        }i\in[\tau]\end{array}
        &\Bigg\}\times 2^Y~.
\end{array}
\end{equation}
where $X = \{\kappa+\tau+1,\ldots,n'\}$ and $Y =
 \{n'+1,\ldots,n\}$.
\end{theorem}
\begin{figure}
\centering \fbox{\includegraphics{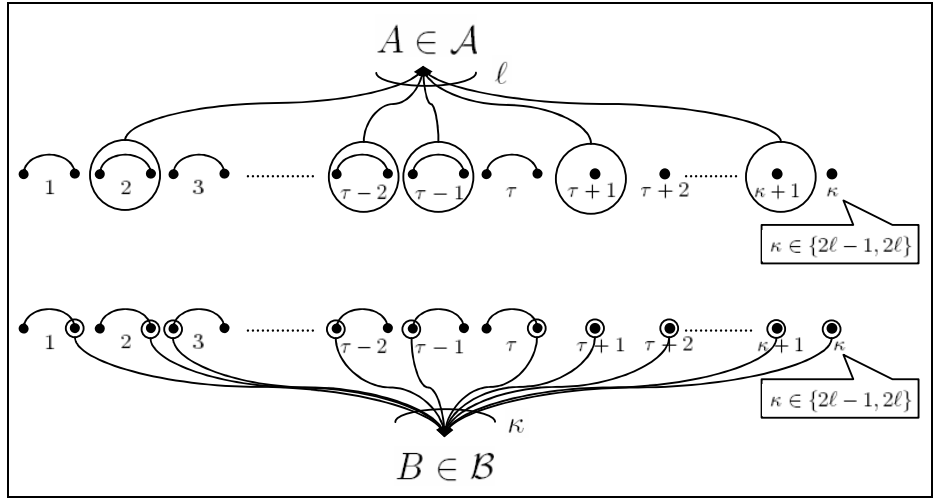}} \caption{The
extremal family \eqref{eq-opt-pair} of $\ell$-cross-intersecting
pairs $\calA,\calB$ in case $n=\kappa+\tau$.}
\label{fig::extremal-family}
\end{figure}
An illustration of the family of extremal pairs $\calA,\calB$
described in Theorem \ref{thm-1} appears in Figure
\ref{fig::extremal-family}. Indeed, this family satisfies:
$$|\calA||\calB|=\binom{\kappa}{\ell}\cdot 2^{|X|} \cdot 2^{\tau+|Y|}=\binom{\kappa}{\ell}2^{n-\kappa}
=\binom{2\ell}{\ell}2^{n-2\ell}~,$$ where the last inequality is by
the choice of $\kappa\in\{2\ell-1,2\ell\}$. The construction of
\cite{ACZ} fits the special case $\tau=0$, $\kappa=2\ell$.

The proof of Theorem \ref{thm-1} combines tools from linear algebra
with techniques from extremal combinatorics, including the
Littlewood-Offord Lemma, extensions of Sperner's Theorem and some
large deviation estimates.

The rest of this paper is organized as follows: Section
\ref{sec::preliminaries} includes some of the ingredients needed for
the proof of Theorem \ref{thm-1}. In order to prove the main result,
we first prove a weaker version of Theorem \ref{thm-1}, which states
that $P_\ell(n) \leq 2^{n+3}/\sqrt{\ell}$ for every sufficiently
large $\ell$ (note that this result alone gives a positive answer to
the above question of Sgall). This is shown in Section
\ref{sec::thm-const}. In Section \ref{sec::thm-1} we reduce the
proof of Theorem \ref{thm-1} to two lemmas, Lemma \ref{lem-rA=o(n)}
and Lemma \ref{lem-rA+sA=Omega(n)}. These lemmas are proved in
Sections \ref{sec::lem-rA=o(n)} and \ref{sec::lem-rA+sA=Omega(n)}
respectively. Section \ref{sec::concluding} contains some concluding
remarks and open problems.

Throughout the paper, all logarithms are in base 2.

\section{Preliminary Sperner-type
Theorems}\label{sec::preliminaries}
\subsection{Sperner's Theorem and the Littlewood-Offord Lemma}
If $P$ is a finite partially ordered set, an antichain of $P$ is a
set of pairwise incomparable elements. Sperner's Theorem
\cite{Sperner} provides a tight upper bound on the maximal size of
an antichain, when $P$ is the collection of all subsets of an
$n$-element set with the subset relation ($A \leq B$ iff $A \subset
B$):
\begin{theorem}[\cite{Sperner}] If $\calA$ is an
antichain of an $n$-element set, then $|\calA| \leq
\binom{n}{\lfloor n/2\rfloor}$.\end{theorem} In
\cite{LittlewoodOfford}, Littlewood and Offord studied a problem
which has the following formulation in the $1$-dimensional case: let
$a_1,\ldots,a_n \in \mathbb{R}$ with $|a_i|
> 1$ for all $i$. What is the maximal number of sub-sums $\sum_{i\in I}a_i$, $I
\subset [n]$, which lie in an interval of length $1$? An immediate
lower bound is $\binom{n}{\lfloor n/2\rfloor}$, when, for some
$\alpha > 1$, half of the $a_i$-s is equal to $\alpha$ and the other
half is equal to $-\alpha$.

Using Sperner's Theorem, Erd\H{o}s \cite{Erdos} gave a tight upper
bound of $\binom{n}{\lfloor n/2\rfloor}$ for the $1$-dimensional
case of the so-called Littlewood-Offord Lemma. To see this, consider
the maximal number of sub-sums of $a_1,\ldots,a_n$, which all belong
to some unit interval. Without loss of generality, we may assume
that all the $a_i$-s are positive (possibly shifting the target unit
interval). Therefore, $a_i > 1$ for all $i$, implying that the
desired family of subsets is an antichain. The result now follows
from Sperner's Theorem. Using a similar argument, Erd\H{o}s proved
the following stronger result:
\begin{lemma}[\cite{Erdos}]\label{lem-littlewood-offord}
Let $a_1,\ldots,a_n \in \mathbb{R}\setminus\{0\}$, and let
$\delta=\min\{|a_i|\}$. Let $T$ be a union of $m$ half-open
intervals, each of width at most $\delta$. Then the number of
sub-sums $\sum_{i\in I}a_i$, $I\subset [n]$, which belong to $T$, is
at most the sum of the $m$ middle binomial coefficients in $n$.
\end{lemma}

\subsection{A bipartite extension of Sperner's Theorem}
The following lemma gives an upper bound on the size of an antichain
of $[n]$, which satisfies an additional requirement with respect to
a pre-defined partition of $[n]$ into into two sets.

\begin{lemma}\label{lem-2-part-antichain} Let $U=[u]$
and $V = [n] \setminus U$, $u \leq n$. If $\calA$ is an antichain of
$[n]$, and in addition satisfies: $|A \cap V| = f(|A \cap U|)$,
where $f:\mathbb{N}\to\mathbb{N}$ is some monotone increasing
function, then $|\calA| \leq \binom{u}{\lfloor
u/2\rfloor}\binom{n-u}{\lfloor (n-u)/2\rfloor}$.
\end{lemma}
The above lemma will follow from the next generalization of
Sperner's Theorem:
\begin{proposition}\label{prop-bipartite-imcomparable-family}
Let $U=[u]$ and $V = [n] \setminus U$, $u \leq n$. If every two sets
$A\neq B\in \calA$ satisfy that either $A\cap U$, $B\cap U$ are
incomparable or $A\cap V$, $B \cap V$ are incomparable, then
$|\calA| \leq \binom{u}{\lfloor u/2\rfloor}\binom{n-u}{\lfloor
(n-u)/2\rfloor}$.
\end{proposition}
\begin{proof} Notice that the upper bound is tight, as it is achieved by
a cartesian product of maximal antichains of $U$ and $V$. The proof
is based on Lubbell's proof \cite{Lubell} of Sperner's Theorem and
the LYM inequality. For each $A \in \calA$, let:
\begin{equation} \label{eq-Au-Av-def} A_U = A \cap U~,~A_V = \{ x-u
: x \in A \cap V\}~.
\end{equation} Let $\sigma \in S_u$
and $\pi \in S_{n-u}$ (where $S_m$ is the symmetric group on $m$
elements) denote two random permutations, chosen uniformly and
independently . We define the event $E_A$ for $A \in \calA$ to be:
$$ E_A = \left( ~A_U=\{\sigma(1),\ldots,\sigma(|A_U|)\} ~\wedge~
A_V = \{\pi(1),\ldots,\pi(|A_V|)\}~\right)~,$$ that is, the first
entries of $\sigma$ form $A_U$, and the first entries of $\pi$ form
$A_V$. The key observation is that the events $E_A$ and $E_B$ are
disjoint for all $A\neq B \in \calA$. To see this, assume that $E_A
\wedge E_B$ holds for some $A \neq B \in \calA$. The fact that the
first entries of $\sigma$ form both $A_U$ and $B_U$ implies that
either $A_U \subset B_U$ or $B_U \subset A_U$, and the same applies
to $A_V,B_V$. Therefore, the assumption on $\calA$ implies that the
events $E_A$ and $E_B$ are indeed disjoint, and thus:
$$ \sum_{A \in \calA} \Pr[E_A] = \Pr[ \bigcup_{A \in \calA} E_A ]
\leq 1~.$$ Since:
$$ \Pr[E_A] = \frac{1}{\binom{u}{|A_U|}\binom{n-u}{|A_V|}}~,$$
it follows that: \begin{equation}\label{eq-gen-LYM-ineq} \sum_{A \in
\calA}\frac{1}{\binom{u}{|A_U|}\binom{n-u}{|A_V|}} \leq
1~.\end{equation}
 Note that in the special case $u=n$ this is the
LYM inequality. The left hand side of \eqref{eq-gen-LYM-ineq} is at
most $\sum_{A \in\calA} 1/\left({\binom{u}{\lfloor
u/2\rfloor}\binom{n-u}{\lfloor (n-u)/2\rfloor}}\right)$ and the
desired result follows.
\end{proof}
\begin{proof}[Proof of Lemma \ref{lem-2-part-antichain}]
Following the notation of Proposition
\ref{prop-bipartite-imcomparable-family}, define $A_U$ and $A_V$ for
each $A\in\calA$ as in \eqref{eq-Au-Av-def}. By Proposition
\ref{prop-bipartite-imcomparable-family}, it suffices to show that,
for all $A\neq B \in \calA$, either $A_U$, $B_U$ are incomparable or
$A_V$, $B_V$ are incomparable. Assume the contrary, and let $A\neq B
\in \calA$ be a counterexample. Without loss of generality, assume
that $A_U \subset B_U$. If $A_V \subset B_V$ then $A \subset B$,
contradicting the fact that $\calA$ is an antichain. It follows that
$B_V \subsetneqq A_V$, and since $f$ is monotone increasing, the
following holds:
$$ |A_V| > |B_V| = f(|B_U|) \geq f(|A_U|)~,$$
contradicting the assumption that $|A_V|=f(|A_U|)$.
\end{proof}

\section{An upper bound tight up to a constant}\label{sec::thm-const}
In this section we prove a weaker version of Theorem \ref{thm-1},
whose arguments will be later extended to prove the precise lower
bound.
\begin{theorem}\label{thm-const}
For any sufficiently large $\ell \in \mathbb{N}$, every
$\ell$-cross-intersecting pair $\calA,\calB \subset
2^{[n]}$ satisfies:
\begin{equation}\label{eq-cross-bound}
|\calA| |\calB| \leq\frac{2^{n+3}}{\sqrt{\ell}}~.
\end{equation}
\end{theorem}

\begin{proof} Let $\calA$ and $\calB$ be as above. A key observation is the
following: it is sufficient to prove \eqref{eq-cross-bound} for the
case where both $\calA$ and $\calB$ are antichains. This follows
from an induction on $n$, where in the case $n=\ell$,
$|\calA||\calB|=1$ and \eqref{eq-cross-bound} clearly holds. Indeed,
suppose that there exist $A_1,A_2\in \calA$ such that $A_1 \subset
A_2$. As $(\calA,\calB)$ are $\ell$-cross-intersecting, this implies
that:
\begin{equation}\label{eq-induction-move}
B \cap (A_2\setminus A_1) = \emptyset \mbox{ for all }B\in\calB~,
\end{equation} hence the restriction of the families $(\calA,\calB)$
to $[n]\setminus (A_2\setminus A_1)$, $(\calA',\calB')$, is an
$\ell$-cross-intersecting pair of an $n'$-element set,
where $n'<n$. By \eqref{eq-induction-move}, $|\calB'|=|\calB|$, and
by the induction hypothesis:
$$|\calA||\calB| \leq 2^{n-n'}|\calA'||\calB'| \leq
\frac{2^{n+3}}{\sqrt{\ell}}~,$$ as required.

For any subset $A\subset [n]$, let $\chi_A \in \{0,1\}^n$ denote its
characteristic vector. Let $\mathcal{F_A}$ and $\mathcal{F_B}$
denote the linear subspaces of $\mathbb{R}^n$ formed by the
characteristic vectors of $\calA$ and $\calB$ respectively:
\begin{equation}
  \label{eq-FA-FB-def}\begin{array}{lll}
  \mathcal{F_A} &=& \mathrm{span}(\{\chi_A : A \in \calA \}) \subset
\mathbb{R}^n~,\\
\mathcal{F_B} &=& \mathrm{span}(\{\chi_B : B \in \calB \}) \subset
\mathbb{R}^n~,\end{array}
\end{equation}
and assume without loss of generality that $\dim(\mathcal{F_A}) \geq
\dim(\mathcal{F_B})$. Choose an arbitrary set $B_1 \in \calB$ and
define:
\begin{equation}
  \label{eq-FB'-def}\begin{array}{l}
   \mathcal{F_B'} = \mathrm{span}(\{\chi_B - \chi_{B_1} : B \in
\calB\})~,\\
k = \dim(\mathcal{F_A})~,~h=\dim(\mathcal{F_B'}) \leq
\dim(\mathcal{F_B})~.\end{array}
\end{equation}
By the definition of $\ell$-cross-intersection, it follows that
$\mathcal{F_A},\mathcal{F_B'}$ are two orthogonal linear subspaces
of $\mathbb{R}^n$, and $k+h \leq n$. Note also that $k \geq h$ by
the assumption on $\dim(\mathcal{F_A})$.

Let $M_\calA$ denote the $k\times n$ row-reduced echelon form
matrix, which is the result of performing Gauss elimination on the
row-vectors $\{ \chi_A : A \in \calA\}$ over $\mathbb{R}$, and let
$M_\calB$ denote the corresponding $h\times n$ matrix for the
vectors $\{ \chi_B - \chi_{B_1} : B \in \calB\}$. As
$\mathrm{rank}M_\calA = k$ and $\mathrm{rank}M_\calB = h$, without
loss of generality we have:
$$ M_\calA = \left(\begin{array}{c|c}I_k\;&\; *\end{array}\right)~,~
M_\calB = \left(\begin{array}{c|c}I_h \;&\; *\end{array}\right)~.$$
where $I_r$ denotes the identity matrix of order $r$ (and the order
of the columns in $M_\calA$ and $M_\calB$ is not necessarily the
same). This implies that any linear combination of the rows of
$M_\calA$ which belongs to $\{0,1\}^n$ has precisely two possible
coefficients for each row: $\{0,1\}$, and in particular, $|\calA|
\leq 2^k$. Similarly, $|\calB| \leq 2^h$ (the two possible
coefficients in the affine combination are now determined by the
vector $\chi_{B_1}$), hence $|\calA||\calB| \leq 2^{k+h} \leq 2^n$,
giving the known upper bound of \cite{FranklRodl}. Observe that if
$k+h \leq n-\log{n}$, we get
$$|\calA||\calB| \leq \frac{2^n}{n} ~,$$ and
\eqref{eq-cross-bound} clearly holds. Therefore, recalling that $k
\geq h$, we may assume that:
\begin{equation}\label{eq-k-initial-bound}
\left\{\begin{array}{rcl} \frac{n}{2}-\frac{1}{2}\log n < &k&\\
n - \log n < &k+h& \leq n\end{array}\right. ~.\end{equation} We
claim that the following statement, which clearly implies
\eqref{eq-cross-bound}, holds:
\begin{equation}
  \label{eq-cross-k-h-bound}
  |\calA||\calB| \leq \frac{2^{k+h+3}}{\sqrt{n}}.
\end{equation}
To show this, we need the next lemma, which will be applied once on
$M_\calA,\calA,k$ and once on $M_\calB,\calB,h$, to conclude that a
constant fraction of the rows of $M_\calA$ and $M_\calB$ have
precisely two non-zero entries, $1$ and $-1$.

\begin{lemma}\label{lem-M-submatrix}
Let $M$ denote a $d\times n$ matrix in row-reduced echelon form: $M
= \left(\begin{array}{c|c}I_d\;&\; *\end{array}\right)$, and let
$\mathcal{D}$ denote an antichain of subsets of $[n]$. Assume that:
\begin{enumerate}
\item The characteristic vectors of $\mathcal{D}$ belong to $w +
\mathrm{span}(M)$, the affine subspace formed by some fixed vector
$w\in\{0,1\}^n$ and the span of the rows of $M$.
\item The antichain $\mathcal{D}$ satisfies $|\mathcal{D}| \geq 8 \cdot 2^d / \sqrt{n}$.
\end{enumerate}
Then there exists a subset of $c$ rows of $M$, $C \subset [d]$,
where $c \geq d - \frac{n}{20} - 10\log n$, such that:
\begin{enumerate}
  \item Every row $i$ of $C$ belongs to $\{0,\pm1\}^n \setminus \{0,1\}^n$.
  \item Every column of the $c\times n$ sub-matrix formed by $C$ contains at most $1$ non-zero entry.
\end{enumerate}
\end{lemma}
\begin{proof}
Our first step is to remove a small portion of the rows of $M$, such
that the remaining rows will have at most one non-zero entry in each
column.
\begin{claim}\label{clm-R-rows}
Let $M,\mathcal{D},w$ satisfy the requirements of Lemma
\ref{lem-M-submatrix}. There exists a set of rows $R \subset [d]$
such that $|R| \leq \frac{n}{25}+ 10\log n$, and each column of $M$
has at most one non-zero value in the remaining $d-|R|$ rows.
\end{claim}
\begin{proof}[Proof of Claim]
 Perform the following process of column-selection on $M$: first,
set $M' = M$. If $M'$ has no column with at least $2$ non-zero
entries, the process ends. Otherwise, perform the following step
(step $j$, for $j\geq 1$): \begin{itemize} \item Let $i_j$ denote
the index of a column of $M'$ with a maximal number of non-zero
entries, $r_j$. \item Let $R_j$ denote the set of rows where the
column $i_j$ is non-zero ($|R_j|=r_j$). \item Replace all these rows
in $M'$ by $0$-rows, and continue the process. \end{itemize} The
result is a sequence of indices, $i_1,\ldots,i_t$ ($t \geq 0$) and a
sequence of sets of rows $R_1, \ldots, R_t$ of sizes $r_1 \geq r_2
\geq \ldots \geq r_t
> 1$, such that the column $i_j$ has $r_j$ non-zero values in the rows $R_j$,
and $R_j \cap R_{j'} = \emptyset $ for all $j \neq j'$. Finally, the
sub-matrix formed by removing the rows $R = \cup_{j=1}^t R_j$ from
$M$ has at most $1$ non-zero entry in every column.

Consider affine combinations (with the affine vector $w$) of the
rows of $M$ which produce a $\{0,1\}^n$-vector. As stated above,
each row of $M$ allows precisely two coefficients in such an affine
combination, as the first $d$ columns of $M$ form the identity
matrix. Clearly, the value of the affine combination at index $i_1$
depends precisely on the $r_1$ coefficients of the rows $R_1$. In
general, if we already chose the coefficients for the rows
$\cup_{j'<j}R_{j'}$, then the value of the affine combination at
index $i_j$ depends only on the choice of the $r_j$ coefficients for
the rows $R_j$.

A simple argument will show that for $1 \leq j \leq t$, at most
$\frac{3}{4}$ of the above $2^{r_j}$ combinations of coefficients
for the rows $R_j$ are indeed valid. To this end, recall the
following simple fact, which corresponds to the Cauchy-Davenport
Theorem when $A,B$ are subsets of $\mathbb{Z}/p\mathbb{Z}$ instead
of $\mathbb{R}$: \begin{equation}\label{eq-simple-c-d}|A+B| \geq
|A|+|B|-1 ~\mbox{ for any two finite nonempty }A,B \subset
\mathbb{R}~,\end{equation} where $A+B = \{a+b: a\in A,~ b \in B\}$.
To see this, simply sort the values of $A$ and $B$ by order of
magnitude, then produce distinct sums by iterating first on $A$,
then on $B$.

Suppose we already chose coefficients for the rows
$\cup_{j'<j}R_{j'}$, and consider the column $i_j$. Select $2$
arbitrary rows $u,v \in R_j$, and fix the choice of coefficients for
the remaining $r_j-2$ rows. We are left with a choice between two
coefficients for $u$, yielding two possible values $a_1,a_2$
contributed by $u$ to the index $i_j$. Similarly, the row $v$
contributes one of two possible values $b_1,b_2$ to the index $i_j$.
Setting $A=\{a_1,a_2\}$ and $B=\{b_1,b_2\}$, the above fact implies
that $|A+B|\geq 3$, hence at least one of the $4$ possible
combinations of $u$ and $v$ gives a non-$\{0,1\}$ value in index
$i_j$ of the resulting affine combination. Therefore, at most
$\frac{3}{4}$ of the $2^{r_j}$ combinations for $R_j$ result in a
$\{0,1\}^n$ vector. We conclude that $|\mathcal{D}| \leq
\left(\frac{3}{4}\right)^t 2^d~,$ and hence $t \leq 2\log n$,
otherwise we would get:
$$ |\mathcal{D}| \leq \frac{2^d}{n^{2\log(4/3)}} <
\frac{2^d}{\sqrt{n}}~,$$ contradicting the assumption on
$|\mathcal{D}|$.

After providing an upper bound on $t$, we wish to bound the term
$\sum_{i=1}^t r_i$. Let $0\leq s\leq t$ denote the maximal index
such that $r_s \geq 6$,
that is: \begin{eqnarray} r_1 \geq r_2 \geq \ldots \geq r_s \geq 6~,\nonumber\\
 6 > r_{s+1} \geq r_{s+2} \geq \ldots \geq r_t >
1~.\nonumber\end{eqnarray} As before, we consider the choice of
coefficients for the rows $R_j$ at step $j$, determining the
$i_j$-th entry of the linear combination. By the Littlewood-Offord
Lemma (Lemma \ref{lem-littlewood-offord}), we conclude that there
are at most $2\binom{r_j}{\lfloor r_j/2 \rfloor} <
\frac{2}{\sqrt{\frac{\pi}{2}r_j}} 2^{r_j}$ possible combinations of
the rows $R_j$ which yield a $\{0,1\}$-value in the $i_j$ column
(note that the inequality $\binom{2x}{x}\leq 2^{2x}/\sqrt{\pi x}$
holds for every integer $x \geq 1$, by the improved approximation
\cite{Robbins} of the error term in Stirling's formula). Applying
this argument to $i_1,\ldots,i_s$, we obtain that:
\begin{equation}\label{eq-D-ri-bound} |\mathcal{D}| \leq 2^d
\prod_{i=1}^s \frac{2\sqrt{2/\pi}}{\sqrt{r_i}}~.
\end{equation}
Observe that every $m$ reals $a_1,\ldots,a_m \geq 2$ satisfy:
$$ \prod_{i=1}^m \frac{1}{a_i} \leq \frac{1}{\sum_{i=1}^m a_i}$$
(this follows by induction on $m$ from the fact that $x y \geq x+y$
for $x,y \geq 2$). Therefore, as $r_i \geq 6 > 2\cdot
(2\sqrt{2/\pi})^2$ for $1 \leq i \leq s$, it follows that:
$$\prod_{i=1}^s \frac{2\sqrt{2/\pi}}{\sqrt{r_i}} \leq
\frac{2\sqrt{2/\pi}}{\sqrt{\sum_{i=1}^s r_i}}~.$$ Combining this
with \eqref{eq-D-ri-bound} we obtain that if $\sum_{i=1}^s r_i
> n/25$, then $|\mathcal{D}| < 8\cdot 2^d/\sqrt{n}$, contradicting
the assumption on $|\mathcal{D}|$. Assume therefore that
$\sum_{i=1}^s r_i \leq n/25$. Altogether, we obtain that $R =
\cup_{j=1}^t R_j$ satisfies:
$$ |R| = \sum_{i=1}^t r_i \leq (\sum_{i=1}^s r_i) + 5(t-s) \leq \frac{n}{25}+ 10\log n~. $$
This completes the proof of the claim.
\end{proof}

It remains to deal with rows which do not belong to $\{0,\pm1\}^n
\setminus \{0,1\}^n$. The next claim provides an upper bound on the
number of such rows in $M$:
\begin{claim}\label{clm-S-rows}
Let $M,\mathcal{D},w$ satisfy the requirements of Lemma
\ref{lem-M-submatrix}, and let $R \subset [d]$ be a set of indices
of rows of $M$ as provided by Claim \ref{clm-R-rows}. Let $S$ denote
the set of indices in $[d] \setminus R$ of rows which do not belong
to $\{0,\pm1\}^n\setminus\{0,1\}^n$. Then $|S| < n/100$.
\end{claim}
\begin{proof}[Proof of Claim]
To prove the claim, fix a linear combination $u$ of the rows $[d]
\setminus S$, and consider all the possible combinations of the rows
of $S$ which can be added to $w'=w+u$ to produce vectors of
$\mathcal{D}$. We will show that the number of these combinations is
at most $2^{s}/\sqrt{\pi s/2}$, where $s = |S|$, and the result will
follow from the assumption on $|\mathcal{D}|$.

Put $S = S_{01} \cup S_{\overline{01}}$, where $S_{01} \subset S$ is
the set of indices of rows in $S$ which are $\{0,1\}^n$ vectors, and
$S_{\overline{01}} = S \setminus S_{01}$. Recall that the first $d$
columns of $M$ form the identity matrix, and that $w\in\{0,1\}^n$,
hence the only two coefficients which can be assigned to the row $i$
to produce $\{0,1\}$ values in the $i$-th column are:
\begin{equation}\label{eq-row-i-coeffs}\left\{\begin{array}
  {ll} \{0,1\}&\mbox{if }w_i=0 \\
   \{0,-1\}&\mbox{if }w_i=1
\end{array}\right.~.
\end{equation}
It will be more convenient to have the coefficients $\{0,1\}$ for
all rows of $S$: to obtain this, subtract each row $i\in S$, whose
coefficients are $\{0,-1\}$, from $w'$, and let $w''$ denote the
resulting vector.

 Let $i\in S_{\overline{01}} $ be an index of a row which does
not belong to $\{0,\pm1\}^n$, and let $j$ denote a column such that
$M_{i j} = \lambda \notin \{0,\pm 1\}$. Crucially, $S \cap R =
\emptyset$, hence column $j$ contains at most one non-zero entry in
the rows of $S$. Therefore, the two possible values of the affine
combination in index $j$ are $\{w'_j,w'_j + \lambda\}$, and as
$0<|\lambda|\neq 1$ it follows that at least one of these values
does not belong to $\{0,1\}$. We deduce that there is at most one
valid choice of coefficients for all the rows $S_{\overline{01}}$.
Denoting this unique combination of the rows of $S_{\overline{01}}$
by $v$, it follows that every linear combination of $S$ which, when
added to $w'$, belongs to $\mathcal{D}$, is the sum of $z = w''+v$
and a linear combination of $S_{01}$.

It remains to set the coefficients of the rows $S_{01}$, and since
each row of $S_{01}$ has $\{0,1\}$ as its coefficients, we are
considering a sum of a subset of the rows of $S_{01}$. Each of these
rows belongs to $\{0,1\}^n$, and in particular, is non-negative: we
claim that the set of possible subsets of $S_{01}$ is therefore an
antichain. To see this, suppose that two distinct subsets
$X,Y\subset S_{01}$, $X \subset Y$, produce (when added to $z$) two
vectors $x,y\in\mathbb{R}^n$ which correspond to sets in
$\mathcal{D}$. The values of $x,y$ at the indices of $S_{01}$ are
determined by the sets $X,Y$ (in fact, these values are equal to
those of the corresponding characteristic vectors), hence $x \neq
y$. Furthermore, as the rows of $S_{01}$ are non-negative, and $X
\subset Y$, we have $x_i \leq y_i$ for all $i \in [n]$. This
contradicts the fact that $\mathcal{D}$ is an antichain. Let
$s'=|S_{01}|$; Sperner's Theorem gives:
$$ |\mathcal{D}| \leq 2^{d-s} \cdot \binom{s'}{\lfloor s'/2\rfloor} \leq
2^{d-s} \cdot \binom{s}{\lfloor s/2\rfloor} \leq
\frac{2^d}{\sqrt{\pi s/2}} ~,$$ and by the assumption on
$|\mathcal{D}|$, we obtain that $s \leq n/100$, completing the proof
of the claim.
\end{proof}

Altogether, Claims \ref{clm-R-rows} and \ref{clm-S-rows} imply that
we can delete at most
$$|R| + |S| \leq \frac{n}{20}+ 10\log n $$ rows of
$M$, and obtain a subset of $c$ rows, $d - \frac{n}{20} - 10\log n
\leq c \leq d$, satisfying the statements of the lemma.
\end{proof}
Note that the requirements of Lemma \ref{lem-M-submatrix} are
satisfied both by $M_\calA,\calA$ and by $M_\calB,\calB$. Indeed, if
either $|\calA| < \frac{8}{\sqrt{n}} \cdot 2^k$ or $|\calB| <
\frac{8}{\sqrt{n}}\cdot 2^h$, then \eqref{eq-cross-k-h-bound} holds
and we are done. The remaining requirement on the characteristic
vectors of $\mathcal{D}$ is satisfied by definition (for $\calA$,
$w$ is the zero vector, whereas for $\calB$, $w = \chi_{B_1}$).

Applying Lemma \ref{lem-M-submatrix} to $M_\calA,\calA$, we obtain a
set of at least $c_1 \geq k - \frac{n}{20}-10\log n$ rows,
$C_1\subset [k]$, such that each row has an entry of $-1$ at some
index $j>k$, and each column has at most $1$ non-zero entry in these
rows. In particular, we get: $ c_1 \leq n - k$, and thus:
$$ k - \frac{n}{20}-10\log n \leq n - k~,$$
and by \eqref{eq-k-initial-bound} we get:
\begin{equation}\label{eq-k-h-bound} \left\{\begin{array}{rcl}
\frac{n}{2} - \log n &\leq k \leq& \frac{21}{40}n + 5 \log n
\\
\frac{19}{40}n-6\log n &\leq h \leq &\frac{n}{2}
\end{array}\right.~.
\end{equation} Next, let $C'_1 \subset C_1$ denote the set of indices of rows
of $C_1$ with precisely two non-zero entries. Notice that, as each
of the columns $\{k+1,\ldots,n\}$ contains at most $1$ non-zero
entry in the rows $C_1$, and on the other hand, each of the rows
$C_1$ contains a non-zero value in one of these columns, it follows
that $|C_1 \setminus C'_1| \leq n-k-c_1$. The lower bound on $c_1$
and \eqref{eq-k-h-bound} give the following bound on $c'_1 =
|C'_1|$:
\begin{equation}\label{eq-c'1-bound}c'_1 \geq c_1 - (n - k - c_1) \geq 3k
-\frac{n}{10}-20\log n - n \geq \frac{2}{5}n -23\log
n~.\end{equation} Since each row $i\in C'_1$ has precisely $2$
non-zero entries, it follows that it has the entry $1$ at index $i$
and the entry $-1$ at some index $j > k$.

Applying Lemma \ref{lem-M-submatrix} to $M_\calB$ and $\calB$, we
obtain a set of at least $c_2 \geq h - \frac{n}{20}-10\log n$ rows,
$C_2 \subset [h]$, and a similar argument to the one above implies
that at most $n-h-c_2$ rows can contain more than $2$ non-zero
entries. Let $C'_2 \subset C_2$ denote the set indices of rows of
$C_2$ with precisely two non-zero entries, and let $c'_2 = |C'_2|$.
By the lower bound on $c_2$ and \eqref{eq-k-h-bound} we obtain:
\begin{equation}\label{eq-c'2-bound}c'_2 \geq c_2 - (n-h-c_2) \geq
3h - \frac{n}{10} - 20\log n - n \geq \frac{13}{40}n - 38\log n~.
\end{equation}
Note that each row $i \in C'_2$ has the entry $1$ at the index $i$
and the entry $-1$ at some index $j > h$.

Finally, notice that \eqref{eq-c'1-bound} and \eqref{eq-c'2-bound}
imply that $c'_1 + c'_2 > n/2$ for a sufficiently large value of
$n$. However, as the rows of $M_\calA$ and $M_\calB$ are orthogonal,
the non-zero entries of each pair of rows $i\in C_1'$ and $j \in
C_2'$ must be in pairwise disjoint columns. In particular, we obtain
that $ 2c'_1 + 2c'_2 \leq n $, yielding a contradiction. Thus,
either $\calA$ or $\calB$ does not meet the requirements of Lemma
\ref{lem-M-submatrix}, and we deduce that \eqref{eq-cross-k-h-bound}
holds.
\end{proof}

\section{Proof of Theorem \ref{thm-1} and two
lemmas}\label{sec::thm-1} Let $\calA$ and $\calB$ denote an
$\ell$-cross-intersection pair of families in $2^{[n]}$. Recall that
in the proof of Theorem \ref{thm-const}, we argued that if, for
instance, $\calA$ is not an antichain, then $\bigcup_{B\in\calB}B
\neq [n]$ (see \eqref{eq-induction-move}). In such a case, letting
$i\in[n]$ be so that $i\notin B$ for all $B\in\calB$, it follows
that $\calA=\calA' \cup \{A \cup \{i\}:A\in\calA'\}$ and
$\calB=\calB'$, where $(\calA',\calB')$ is an optimal
$\ell$-cross-intersecting pair on $[n]\setminus\{i\}$. Therefore, by
induction, the structure of $\calA,\calB$ is as specified in Theorem
\ref{thm-1}, where the parameter $n'$ (determining the set $X$ in
\eqref{eq-opt-pair}) accounts for the modification of
$(\calA',\calB')$ to $(\calA,\calB)$. The same consideration applies
when $\bigcup_{A\in\calA}A \neq [n]$, which follows when $\calB$ is
not an antichain (in this case, the set $Y$ in \eqref{eq-opt-pair}
treats the modification of $\calB'$ to $\calB$). Altogether, we may
assume that $\calA,\calB$ are both antichains, and furthermore:
\begin{equation}
  \label{eq-A,B-contain-all-elements}
  \bigcup_{A\in\calA}A =   \bigcup_{B\in\calB}B = [n]~.
\end{equation}
It remains to prove that in this case $|\calA||\calB| \leq
\binom{2\ell}{\ell}2^{n-2\ell}$, and that equality holds iff for
some
\begin{equation}\label{eq-opt-pair-params'}
\begin{array}
  {l}
  \kappa\in\{2\ell-1,2\ell\}~,~\tau \in \{0,\ldots,\kappa\}~,\\
  \kappa+\tau = n,
\end{array}\end{equation}
the following holds up to a relabeling of the elements of $[n]$ and
swapping $\calA,\calB$:
\begin{equation} \label{eq-opt-pair'}
\begin{array}{lll}
        \mathcal{A} &= \Bigg\{
    \displaystyle{\bigcup_{T \in J} T} ~:~J
    \subset\bigg\{\begin{array}{c}
           \{1,\kappa+1\},\ldots,\{\tau,\kappa+\tau\},\\
           \{\tau+1\},\ldots,\{\kappa\}
    \end{array}\bigg\}~,~
    |J|=\ell
    &\Bigg\}~,\\
    \\
        \mathcal{B}&=\Bigg\{L \cup \{\tau+1,\ldots,\kappa\} :
        \begin{array}{l}
        L \subset \{1,\ldots,\tau,\kappa+1,\ldots,\kappa+\tau\}\\
        |L\cap \{i,\kappa+i\}| = 1 \mbox{ for all
        }i\in[\tau]\end{array}
        &\Bigg\}~.
\end{array}
\end{equation}

Following the notations of Theorem \ref{thm-const}, define
$\mathcal{F_A},\mathcal{F_B'},k,h$ as in \eqref{eq-FA-FB-def} and
\eqref{eq-FB'-def}, obtaining $k \geq h$. Recall that the proof of
Theorem \ref{thm-const} implies that $|\calA||\calB| \leq
2^{k+h+3}/\sqrt{n}$ provided that $\ell$ is sufficiently large
(equation \eqref{eq-cross-k-h-bound}). This implies that if $k+h
\leq n - 4$ then:
$$ |\calA||\calB| \leq \frac{1}{2}\cdot \frac{2^n}{\sqrt{n}}~,$$
and as $\frac{1}{2} < 1/\sqrt{\pi}$, the pair $\calA,\calB$ is
suboptimal. Assume therefore that $k + h \geq n-3$:
\begin{equation}\label{eq-k-h-revised-bound}
\left\{\begin{array}{rcc} \frac{n-3}{2} \leq &k&\\
n - 3 \leq &k+h& \leq n\end{array}\right. ~.\end{equation} Observe
that, as the rows of $M_\calA$ are orthogonal to the rows of
$M_\calB$, we may assume without loss of generality that:
$$ M_\calA = \left(\begin{array}{c|c}I_k\;&\; *\end{array}\right)~,~
M_\calB = \left(\begin{array}{c|c}* \;&\; I_h\end{array}\right)~.$$
To see this, first perform Gauss elimination on a basis for
$\mathcal{F_A}$ to obtain $M_\calA$. Next, perform Gauss elimination
on a basis for $\mathcal{F_B'}$, and notice that, as the rows of
$M_\calA$ and $M_\calB$ are pairwise orthogonal, it is always
possible to find a leading non-zero entry at some index $j > k$.
Once $M_\calB$ is in row-reduced echelon form, we may relabel the
elements $k+1,\ldots,n$ to obtain the above structure.

 We again apply the arguments of Lemma \ref{lem-M-submatrix} on
$\calA,M_\calA$ and on $\calB,M_\calB$, only this time we perform
the calculations more carefully. Let $R_A \subset [k]$ denote the
subset of the rows of $M_\calA$ which are selected by the process
described in Claim \ref{clm-R-rows}. That is, we repeatedly select
an arbitrary column with at least $2$ non-zero entries, while one
exists, add the rows where it is non-zero to $R_A$, and delete them
from $M_\calA$. While in Claim \ref{clm-R-rows} we repeatedly
selected a column with a maximal number of non-zero entries, here we
allow an arbitrary choice when selecting the next column with at
least $2$ non-zero entries. Let $r_A = |R_A|$, and define $R_B
\subset [h]$ and $r_B = |R_B|$ similarly for $M_\calB$.

Let $S_A \subset [k]\setminus R_A$ denote the indices of rows of
$M_\calA$, which belong neither to $R_A$ nor to
$\{0,\pm1\}^n\setminus \{0,1\}^n$. That is, $S_A$ denotes the rows
which were treated by Claim \ref{clm-S-rows}. Let $s_A = |S_A|$, and
define $S_B \subset [h] \setminus R_B$ and $s_B = |S_B|$ similarly
for $M_\calB$.

The following lemma, proved in Section \ref{sec::lem-rA=o(n)},
determines the optimal pairs $\calA,\calB$ when $r_A+s_A = o(n)$:
\begin{lemma}\label{lem-rA=o(n)} If there exists some order of column
selection when producing the set $R_A$ such that $r_A + s_A = o(n)$,
then $|\calA||\calB| \leq \binom{2\ell}{\ell}2^{n-2\ell}$.
Furthermore, equality holds iff either:
\begin{equation}
  \label{eq-rA-sA=o(n)-opt-family-1}
\begin{array}{l}
 \mbox{\small$M_\calA = \left( \begin{array}{c|c|c||c}
  I_{k-1} & \begin{array}{c}
    0\\ \vdots \\ 0
  \end{array} & -I_{k-1} & \begin{array}{c}
    0\\ \vdots \\ 0
  \end{array} \\
  \hline
  0 & 1 & 1 \ldots 1 & 1
  \end{array}\right) ~,~
    M_\calB = \left( \begin{array}{c|c|c||c}
  I_{k-1} & \begin{array}{r}
    -1\\ \vdots \\ -1
  \end{array} & I_{k-1} & \begin{array}{c}
    0\\ \vdots \\ 0
  \end{array} \\
  \hline\hline
  0 & -1 & 0\ldots 0 & 1
  \end{array}\right)$}\\
  h \in \{2\ell-2,2\ell-1\}~,~h+k=n~,~k \in \{\frac{n}{2},\frac{n+1}{2}\}~,~
  B_1 = \cup_{i\in[\ell]} \{(i,k+i)\}
    \end{array}
  \end{equation}
or :
\begin{equation}
  \label{eq-rA-sA=o(n)-opt-family-2}
 \begin{array}{l}
 \mbox{\small$M_\calA = \left( \begin{array}{c|c|c|c||c|c}
  I_{k-2} & \begin{array}{c}
    0\\ \vdots \\ 0
  \end{array} & \begin{array}{c}
    0\\ \vdots \\ 0
  \end{array} & -I_{k-2} & \begin{array}{c}
    0\\ \vdots \\ 0
  \end{array} & \begin{array}{c}
    0\\ \vdots \\ 0
  \end{array}\\
  \hline
  0 & 1 & 0 & 1 \ldots 1 & 1 & 1 \\
  0 & 0 & 1 & 1 \ldots 1 & 1 & 1
  \end{array}\right) ,
    M_\calB = \left( \begin{array}{c|c|c|c||c|c}
  I_{k-2} & \begin{array}{r}
    -1\\ \vdots \\ -1
  \end{array} & \begin{array}{r}
    -1\\ \vdots \\ -1
  \end{array} & I_{k-2} & \begin{array}{c}
    0\\ \vdots \\ 0
  \end{array} & \begin{array}{c}
    0\\ \vdots \\ 0
  \end{array} \\
  \hline \hline
  0 & -1 & -1 & 0\ldots 0 & 1 & 0\\
    0 & -1 & -1 & 0\ldots 0 & 0 & 1\\
  \end{array}\right)$}\\
  h \in \{2\ell-2,2\ell-1\}~,~h+k=n~,~k \in \{\frac{n}{2},\frac{n+1}{2},\frac{n}{2}+1\}~,~
  B_1 = \cup_{i\in[\ell]} \{(i,k+i)\}
    \end{array}
  \end{equation}
up to a relabeling of the elements of $[n]$ and the choice of $B_1$.
In both cases above, the pair $(\calA,\calB)$ belongs to the family
\eqref{eq-opt-pair'} with $\kappa=h+1$, $\tau = k-1$ and swapping
$\calA,\calB$.
\end{lemma}
In the above figures \eqref{eq-rA-sA=o(n)-opt-family-1} and
\eqref{eq-rA-sA=o(n)-opt-family-2}, the columns to the right of the
double-line-separators and the rows below the double-line-separators
appear or not, depending on the value of $k$.

The remaining case is treated by the next lemma, which is proved in
Section \ref{sec::lem-rA+sA=Omega(n)}, and concludes the proof of
the theorem:
\begin{lemma}
  \label{lem-rA+sA=Omega(n)}
  If every order of column selection when producing the set $R_A$
  gives $r_A + s_A = \Omega(n)$, then $|\calA||\calB| \leq
  \binom{2\ell}{\ell}2^{n-2\ell}$. Furthermore, equality holds iff:
  \begin{equation}
  \label{eq-rA-sA=Omega(n)-opt-family}
 \begin{array}{l}
 M_\calA = \left( \begin{array}{c|c|c}
  I_h & 0 & I_h \\
  \hline
  0 & I_{k-h} & 0 \end{array}\right) ~,~
    M_\calB = \left( \begin{array}{c|c|c}
    -I_h & 0 & I_h\end{array}\right)\\
  k \in \{2\ell-1,2\ell\}~,~h+k=n~,~
  B_1 = [k]
    \end{array}
  \end{equation}
up to a relabeling of the elements of $[n]$. In this case, the pair
$(\calA,\calB)$ belongs to the family \eqref{eq-opt-pair'} with
$\kappa= k $ and  $\tau = h$.
\end{lemma}

\begin{remark}
  It is, in fact, not difficult to check that if, in one order of
  column selection we have $r_A + s_A = \Omega(n)$, so is the case in any
  order, but the above formulation suffices for our purpose.
\end{remark}
\section{Proof of Lemma \ref{lem-rA=o(n)}}\label{sec::lem-rA=o(n)}
Let $C_1 = [k] \setminus (R_A \cup S_A)$. By the assumption on $r_A,
s_A$ and the fact that $k\geq \frac{n-3}{2}$ we deduce that $|C_1| =
(1-o(1))k$. Recall that each column of $M_\calA$ contains at most
one non-zero entry in the rows of $C_1$, and that each row of $C_1$
belongs to $\{0,\pm1\}^n\setminus \{0,1\}^n$. Hence, $n \geq k +
|C_1| = (2-o(1))k$. Altogether, we obtain that:
\begin{equation}
  \label{eq-k-h-n/2}
  k = \left(\frac{1}{2}+o(1)\right)n~~,~~h = \left(\frac{1}{2}-o(1)\right)n~.
\end{equation}
The $\{1,-1\}$ entries in each row of $C_1$ account for $2|C_1| =
(1-o(1))n$ distinct columns, leaving at most $o(n)$ columns which
may contribute additional values to rows of $C_1$. Again, as each
column contains at most $1$ non-zero entry in the rows of $C_1$, the
set of all rows with non-zero entries either in these columns, or in
columns $\{k+1,\ldots,n-h\}$ (at most $3$ columns), is of size
$o(n)$. We obtain that, without loss of generality:
\begin{equation}\label{eq-n/2-A-structure}
\begin{array}{c|c|c|c|c}
\multicolumn{2}{r@{\mbox{\tiny$ \dashrightarrow|$}}}{\mbox{\tiny$
\dashleftarrow \cdots \quad k \quad \cdots $}} &
\multicolumn{1}{@{\mbox{\tiny$|$}}c@{\mbox{\tiny$|$}}}{\mbox{\tiny$\dashleftarrow
\leq3\dashrightarrow$}} &
\multicolumn{2}{@{\mbox{\tiny$|\dashleftarrow$}}l}{\mbox{\tiny$\cdots h \cdots\dashrightarrow$}} \\
M_\calA = \left(\begin{array}{c}I_{k'}\\ 0 \end{array}\right. &
\begin{array}{c}0\\ I_{k-k'}
\end{array} &
\begin{array}{c}0\\ * \end{array} &
\begin{array}{c}-I_{k'}\\
*
\end{array} &
\left.\begin{array}{c}0\\ * \end{array}
\right)~,\end{array}\end{equation}
 where $k' = (1-o(1))k=(1-o(1))h$.
The above structure of $M_\calA$ provides a quick bound on
$|\calA|$. Consider column $n-h+1$; if this column contains at least
$2$ non-zero entries, then we gain a factor of $\frac{3}{4}$ by
\eqref{eq-simple-c-d}. Otherwise, the fact that $M_{n-h+1,1}=-1$
implies that the coefficient of row $1$ is necessarily $0$, giving a
factor of $\frac{1}{2}$. Therefore:
\begin{equation}
  \label{eq-rA-sA=o(n)-calA-3/4-bound}
  |\calA| \leq \frac{3}{4} \cdot 2^k~.
\end{equation}
For another corollary of \eqref{eq-n/2-A-structure}, notice that for
all $i\in[k']$, row $i$ of $M_\calA$ contains $1,-1$ in columns
$i,n-h+i$ respectively (and $0$ in the remaining columns), and is
orthogonal to all rows of $M_\calB$. It follows that columns
$i,n-h+i$ are equal in $M_\calB$ for all $i \in [k']$, and hence:
\begin{equation}\label{eq-n/2-B-structure}
\begin{array}{c|c|c|c|c}
\multicolumn{2}{r@{\mbox{\tiny$ \dashrightarrow|$}}}{\mbox{\tiny$
\dashleftarrow \cdots k \cdots $}} &
\multicolumn{1}{@{\mbox{\tiny$|$}}c@{\mbox{\tiny$|$}}}{\mbox{\tiny$\dashleftarrow
\leq3\dashrightarrow$}} &
\multicolumn{2}{@{\mbox{\tiny$|\dashleftarrow$}}l}{\mbox{\tiny$\cdots \quad h \quad \cdots\dashrightarrow$}} \\
M_\calB = \left(\begin{array}{c}I_{k'}\\ 0 \end{array}\right. &
\begin{array}{c}*\\ *
\end{array} &
\begin{array}{c}*\\ * \end{array} &
\begin{array}{c}I_{k'}\\
0
\end{array} &
\left.\begin{array}{c}0\\ I_{h-k'} \end{array}
\right)~.\end{array}\end{equation}
 We claim that the above structure
of $M_\calB$ implies that $r_B + s_B = (1-o(1))h$. Indeed, once we
delete the rows $R_B \cup S_B$ from $M_\calB$, each row must contain
an entry of $-1$, which must reside in one of the columns
$k'+1,\ldots,n-h$. As each column contains at most one non-zero
entry in rows $[h]\setminus (R_B \cup S_B)$, we deduce that $n-h-k'
\geq h - r_B - s_B$, and equivalently: $$r_B + s_B \geq 2h + k' - n
= (1-o(1))h = \left(\frac{1}{2}-o(1)\right)n~,$$ where the last two
equalities are by \eqref{eq-k-h-n/2} and the fact that
$k'=(1-o(1))k$. Recall that the analysis of Claim \ref{clm-R-rows}
implies that, if $R_B$ is nonempty, then at most $2^{r_B+1}
/\sqrt{\frac{\pi}{2} r_B}$ linear combinations of the rows of $R_B$
are valid in order to produce a $\{0,1\}^n$ vector from the rows of
$M_\calB$. Furthermore, if $S_B$ is nonempty, then for each choice
of coefficients for the rows $[h]\setminus S_B$, Claim
\ref{clm-S-rows} implies that at most $2^{s_B} /\sqrt{\frac{\pi}{2}
s_B}$ combinations of the rows of $S_B$ are valid in order to
produce a $\{0,1\}^n$ antichain of vectors from the rows of
$M_\calB$. Since in our case we have $r_B + s_B = \Omega(n)$, at
least one of $r_B,s_B$ is $\Omega(n)$, and we deduce that:
\begin{equation}\label{eq-B-sqrt-n-factor}|\calB| = O(2^h / \sqrt{n})~.\end{equation}
Furthermore, if both $r_B = \omega(1)$ and $s_B = \omega(1)$ we get
$|\calB| = O(\frac{2^h}{\sqrt{r_B s_B}}) = o(2^h / \sqrt{n})$ and
hence (regardless of the structure of $M_\calA$) $|\calA||\calB| =
o(2^{k+h} / \sqrt{n}) \leq o(2^n / \sqrt{\ell})$, showing this
cannot be an optimal configuration, as required. The same
consequence is obtained if either $r_A = \omega(1)$ or $s_A =
\omega(1)$, as in this case $|\calA| = o(2^k)$. Assume therefore
that $r_A + s_A = O(1)$, and by the above arguments we obtain that:
\begin{eqnarray}&k = \frac{n}{2}
+ O(1)~,~h = \frac{n}{2} - O(1)~,\\
&k' = k - O(1)~,\\
&r_B = O(1)~,~s_B = h-O(1) ~~\mbox{ or }~~r_B = h-O(1)~,~s_B =
O(1)\label{eq-rB+sB}~.\end{eqnarray} At this point, we claim that
either $n = (4+o(1))\ell$, or the pair $\calA,\calB$ is suboptimal:
\begin{claim}\label{clm-n-geq-4ell} Let $\calA,\calB$ be as above, then either
$|\calA||\calB| = o(2^n / \sqrt{n})$ or $n = (4+o(1))\ell$.
\end{claim}
\begin{proof}
Fix a choice of coefficients for the last $k-k'$ rows of $M_\calA$,
yielding a linear combination $w_A$. By the structure of $M_\calA$
specified in \eqref{eq-n/2-A-structure}, if for some index
$i\in[k']$, $w_A$ does not equal $0$ at index $i$ or does not equal
$1$ at index $n-h+i$, then the $i$-th row of $M_\calA$ has at most
one valid coefficient. Thus, if there are $\omega(1)$ such indices,
we deduce that there are at most
$o(2^{k'})$ combinations of the rows $[k']$ of $M_\calA$ which
extend $w_A$ to an element of $\calA$. Therefore, by
\eqref{eq-B-sqrt-n-factor}, this choice of $w_A$ counts for at most
$o(2^{k'+h}/\sqrt{n})$ pairs $(A,B)\in \calA \times \calB$. Summing
over all $2^{k-k'}$ choices for $w_A$, this amounts to at most
$o(2^n/\sqrt{n})$ pairs $(A,B)\in \calA \times \calB$, and we may
thus assume that at least $k' - O(1)$ of the indices $j \in [k']$
satisfy
\begin{equation}\label{eq-w-A-j-relation}
w_A^{(j)}=0~,~w_A^{(n-h+j)}=1~.\end{equation} Next, fix a choice of
coefficients for the last $h-k'$ rows of $M_\calB$, yielding an
affine combination (together with $\chi_{B_1}$) $w_B$, and consider
the structure of $M_\calB$ specified in \eqref{eq-n/2-B-structure}.
Every index $j \in [k']$ for which $\chi_{B_1}^{(j)} \neq
\chi_{B_1}^{(n-h+j)}$ implies that the row $j$ has at most one valid
coefficient. Thus, if there are $\omega(1)$ such indices, it follows
that $w_B$ can be extended to at most $o(2^h / \sqrt{n})$ elements
of $\calB$. To see this, take $m=\omega(1)$ and yet $m=o(n)$ such
rows, arbitrarily; there is at most one legal combination for these
rows. As $r_B+s_B=\Omega(n)$, the remaining rows have at most
$O(2^{h-m}/\sqrt{n})$ combinations, and the result follows.

Altogether, we may assume that $k'-O(1)$ of the indices $j \in [k']$
satisfy:
\begin{equation}
  \label{eq-chi-B1-j-relation}
  \chi_{B_1}^{(j)} = \chi_{B_1}^{(n-h+j)}~.
\end{equation}
Let $L \subset [k']$ denote the indices of $[k']$ which satisfy both
\eqref{eq-w-A-j-relation} and \eqref{eq-chi-B1-j-relation}. It
follows that $|L| = h - O(1)$, and for each $i \in L$, the choice of
a coefficient for row $i$ exclusively determines between the cases
$i,n+h-i \in B$ and $i,n+h-i \notin B$.

Fix a choice of coefficients for the remaining rows of $M_\calA$,
and let $A$ denote the resulting set, and fix a choice of
coefficients for all rows of $M_\calB$ except those whose indices
are in $L$. For each $i \in L$, let $X_i$ denote the variable whose
value is $1$ if we choose a coefficient for the row $i$ such that
$i,n+h-i \in B$ and $0$ otherwise. Recall that $A$ contains
precisely one element from each pair $\{i,n+h-i~:~i\in L\}$.
Therefore, any choice of coefficients of the rows $L$ in $M_\calB$
gives a set $B$ which satisfies:
\begin{equation}\label{eq-ell-binomial-sum}\ell = |A \cap B| = (\sum_{i
\in L} X_i ) + O(1)~,\end{equation} where the $O(1)$-term accounts
for the intersection of $A$ with at most $n-2|L|=O(1)$ indices.
Choose one of each pair of coefficients for each row of $L$
uniformly at random and independently of the other rows, to obtain
that $X = \sum_{i\in L}X_i$ has a binomial distribution
$\mathrm{Bin}(\frac{n}{2}-O(1),\frac{1}{2})$. Fix some small
$\epsilon
> 0$; by the Chernoff bound (see, e.g., \cite{ProbMethod}, Chapter A.1):
$$\Pr[|X - \frac{n}{4}| > \epsilon n] \leq O\left(\exp(-\Omega(n))\right)~,$$
thus if $|\ell - \frac{n}{4}| > \epsilon n$ then at most $ O(2^h /
\exp(\Omega(n)))$ sets $B \in\calB$ can be produced from $w_B$ and
we are done. We conclude that $\ell = (\frac{1}{4} + o(1))n$.
\end{proof}

The last claim, along with \eqref{eq-rB+sB}, implies that the case
$s_B = h-O(1)$ is suboptimal. Indeed, in this case: $$  |\calB| \leq
\frac{2^h}{\sqrt{\pi s_B / 2}} = (1+o(1))\frac{2^h}{\sqrt{\pi h /2
}} = (1+o(1))\frac{2^h}{\sqrt{\pi n /4 }}
  = (1+o(1))\frac{2^h}{\sqrt{\pi \ell }}~,$$ where the last inequality is by
Claim \ref{clm-n-geq-4ell}. Combining this with
\eqref{eq-rA-sA=o(n)-calA-3/4-bound}, we deduce that
$|\calA||\calB|$ is at most $(\frac{3}{4}+o(1))2^n/\sqrt{\pi\ell}$,
and that the pair $\calA,\calB$ is suboptimal.

It remains to deal with the case $r_B = h - O(1)$, in which case we
have:
\begin{equation}\label{eq-n/2-b-bound}
  |\calB| \leq \frac{2^{h+1}}{\sqrt{\pi r_B /2 }} = (2+o(1))\frac{2^h}{\sqrt{\pi
  \ell}}~,
\end{equation}
and hence ($|\calA|\leq \frac{3}{4}\cdot 2^k$), $|\calA||\calB| \leq
(\frac{3}{2}+o(1))2^{k+h}/\sqrt{\pi\ell}$. If $k+h < n$, it follows
that $|\calA||\calB|$ is at most
$(\frac{3}{4}+o(1))2^n/\sqrt{\pi\ell}$, and again the pair
$\calA,\calB$ is suboptimal. We may thus assume:
$$ k + h = n~,~r_B = h - O(1)~,~s_B = O(1)~.$$
To complete the proof of the lemma, we show that either
$|\calA||\calB| \leq (\delta+o(1))2^n/\sqrt{\pi\ell}$ for some fixed
$\delta < 1$, or all columns of $M_\calB$ except either $1$ or $2$
have at most $1$ non-zero entry, whereas the remaining columns are
of the form $(-1, \ldots,-1)$. This will imply that either
\eqref{eq-rA-sA=o(n)-opt-family-1} holds or
\eqref{eq-rA-sA=o(n)-opt-family-2} holds. For this purpose, we must
first concentrate on the $(k-k')\times k'$ sub-matrix of $M_\calA$,
on rows $\{k'+1,\ldots,k\}$ and columns \{$k+1,\ldots,k+k'$\}. This
sub-matrix appears boxed in diagram
\eqref{eq-n/2-MA-MB-structure-k+h=n}, which reflects the form of
$M_\calA$ and $M_\calB$ given the fact $k+h=n$:
\begin{equation} \label{eq-n/2-MA-MB-structure-k+h=n}
\begin{array}{c|c|c|c}
\multicolumn{2}{r@{\mbox{\tiny$\dashrightarrow|$}}}{\mbox{\tiny$\dashleftarrow
\cdots\cdots k \cdots\cdots $}} &
\multicolumn{2}{@{\mbox{\tiny$|\dashleftarrow$}}l}{\mbox{\tiny$\cdots\cdots h \cdots\cdots \dashrightarrow$}} \\
 M_\calA = \left(\begin{array}{c}I_{k'}\\ 0
\end{array} \right.&
\begin{array}{c}0 \\ I_{k-k'} \end{array} &
\begin{array}{c}-I_{k'}\\ \boxed{~*~} \end{array} &
\left.\begin{array}{c}0\\ * \end{array} \right) \\
\noalign{\medskip}
\multicolumn{2}{r@{\mbox{\tiny$\dashrightarrow|$}}}{\mbox{\tiny$\dashleftarrow
\cdots\cdots k \cdots\cdots $}} &
\multicolumn{2}{@{\mbox{\tiny$|\dashleftarrow$}}l}{\mbox{\tiny$\cdots\cdots h \cdots\cdots \dashrightarrow$}} \\
 M_\calB = \left( \begin{array}{c}I_{k'}\\ 0 \end{array} \right. &
\begin{array}{c}* \\ * \end{array} &
\begin{array}{c}I_{k'}\\ 0 \end{array} &
\left.\begin{array}{c}0\\ I_{h-k'} \end{array} \right)
\end{array}~.
\end{equation}
Suppose the linear combination of rows $k'+1,\ldots,k$ of $M_\calA$
is some vector $w_A$. A key observation is the following: if $w_A$
has $\omega(1)$ entries not equal to $1$ in indices
$\{k+1,\ldots,k+k'\}$, then at most $o(2^{k'})$ combinations of the
remaining rows can be added to $w_A$ to produce a vector in
$\{0,1\}^n$. This follows directly from the structure of $M_\calA$
in \eqref{eq-n/2-A-structure}, as the fact that $w_A^{(k+j)}\neq 1$
forces the coefficient of row $j$ to be $0$. Using the above
observation, we will show that either $|\calA| \leq
(\frac{3}{8}+o(1))2^k$, or at most $O(1)$ columns of $M_\calA$ with
indices $\{k+1,\ldots,k+k'\}$ are not of one of the forms
$\{(-1,1,0,\ldots,0),(-1,1,1,0,\ldots,0)\}$ (at some coordinate
order). Consider the following three cases:
\begin{enumerate}[(I)]
\item \textbf{$\omega(1)$ columns of $M_\calA$ contain at least $3$ non-zero entries
in rows $\{k'+1,\ldots,k\}$:}  Let $S$ denote the indices of columns
in $\{k+1,\ldots,k+k'\}$ for which $M_\calA$ has non-zero entries in
rows $\{k'+1,\ldots,k\}$. The Littlewood-Offord Lemma implies that,
whenever there are $t$ non-zero entries in a single column in these
rows, then at most $m=2^{k-k'-t} \binom{t}{\lfloor t/2\rfloor }$ of
the $2^{k-k'}$ possible linear combinations of these rows can
produce a value of $1$. Notice that for $t\geq 3$ we get
$\binom{t}{\lfloor t/2\rfloor}/2^t \leq \frac{3}{8}$, hence
$m/2^{k-k'}\leq \frac{3}{8}$. Next, let each column which has at
least $3$ non-zero entries in rows $\{k'+1,\ldots,k\}$ ``rate'' $m$
linear combinations, including all those for which it gives a value
of $1$. It follows that choosing any combination for rows
$\{k'+1,\ldots,k\}$ excluding the most popular set of $m$ linear
combinations, yields values not equal to $1$ in at least
$|S|/\binom{2^{k'-k}}{m}=\Omega(|S|)=\omega(1)$ columns, hence (by
the above observation) such combinations contribute $o(2^{k})$
vectors to $\calA$. We deduce that $|\calA| \leq (\frac{3}{8}+o(1))
2^k$.

\item \textbf{$\omega(1)$ columns of $M_\calA$ contain $2$ non-zero
entries $\neq(1,1)$ in rows $\{k'+1,\ldots,k\}$:} The argument here
is similar to the argument in the previous item. If a column has two
non-zero entries $(x,y)\neq(1,1)$ in rows $k'+1,\ldots,k$, then the
possible values of the linear combination at this column are
$\{0,x,y,x+y\}$. At most $1$ of these $4$ values can be $1$, hence
at most $m=2^{k-k'-2}$ of the combinations yield a value of $1$ at
this column. By the above argument, we deduce that $|\calA| \leq
(\frac{1}{4}+o(1)) 2^k$.

\item \textbf{$\omega(1)$ columns of $M_\calA$ contain at most $1$ non-zero
entry $\neq 1$ in rows $\{k'+1,\ldots,k\}$:} this case is the
simplest, following directly from the observation. Indeed, every
linear combination of the rows $k'+1,\ldots,k$ has $\omega(1)$
entries which do not equal $1$ in columns $\{k+1,\ldots,k+k'\}$,
hence $|\calA| = o(2^k)$.
\end{enumerate}
Note that if $|\calA| \leq (\frac{3}{8}+o(1))2^k$, then
$|\calA||\calB|\leq (\frac{3}{4}+o(1))2^n/\sqrt{\pi\ell}$ by
\eqref{eq-n/2-b-bound}, as required. Assume therefore that $M_\calA$
has at most $O(1)$ columns among $\{k+1,\ldots,k+k'\}$, whose set of
non-zero entries in rows $\{k'+1,\ldots,k\}$ is neither $\{1\}$ nor
$\{1,1\}$. We use the abbreviation $\{1\}$-columns and
$\{1,1\}$-columns for the $k'-O(1)$ remaining columns whose non-zero
entries in rows $\{k'+1,\ldots,k\}$ of $M_\calA$ are $\{1\}$ and
$\{1,1\}$ respectively; according to this formulation:
\begin{equation}
  \label{eq-(1)-columns-and-(1,1)-columns}
  k'-O(1)\mbox{ of columns }\{k+1,\ldots,k'\}\mbox{ of $M_\calA$ are $\{1\}$-columns or $\{1,1\}$-columns}~.
\end{equation}
The two cases of whether there are $\omega(1)$ or $O(1)$
$\{1\}$-columns, are treated by Claims \ref{clm-omega(1)-1-columns}
and \ref{clm-O(1)-1-columns} respectively, and determine which of
the two optimal families, stated in
\eqref{eq-rA-sA=o(n)-opt-family-1},\eqref{eq-rA-sA=o(n)-opt-family-2},
is obtained. These two claims are stated and proved in Subsections
\ref{sec::clm-omega(1)-1-columns} and \ref{sec::clm-O(1)-1-columns}.
\subsection{The optimal family
\eqref{eq-rA-sA=o(n)-opt-family-1}}\label{sec::clm-omega(1)-1-columns}
\begin{claim}\label{clm-omega(1)-1-columns}
  If $\omega(1)$ of columns $\{k+1,\ldots,k+k'\}$ of $M_\calA$ are
  $\{1\}$-columns, then \eqref{eq-rA-sA=o(n)-opt-family-1} holds.
\end{claim} \begin{proof} It follows that some row of $\{k'+1,\ldots,k\}$
contains a value of $1$, which is the single non-zero entry of this
column in these rows, in $\omega(1)$ columns of
$\{k+1,\ldots,k+k'\}$ (take the most popular row of
$\{k'+1,\ldots,k\}$). Without loss of generality, assume that this
row is row $k$, the last row of $M_\calA$. By the observation above,
the coefficient for row $k$ of $M_\calA$ must be $1$, otherwise only
$o(2^k)$ combinations of the remaining rows produce vectors in
$\{0,1\}^n$. This has several consequences:
\begin{enumerate}[(1)]
  \item Row $k$ contains the value $1$ in columns
  $\{k+1,\ldots,k+k'\}$. To see this, notice that if $(M_\calA)_{k,k+j} \neq
  1$ for some $j \in [k']$,
then $|\calA|\leq
 \left(\frac{1}{4}+o(1)\right)2^k$: either the coefficient for row
$k$ is $0$, contributing $o(2^k)$ vectors to $|\calA|$, or it is
$1$, forcing the coefficient of row $j$ to be $0$.
 \item Row $k$ contains $\{0,1\}$ values in columns
 $\{k+k'+1,\ldots,n\}$. Indeed, if
$(M_\calA)_{k,k+j}\notin\{0,1\}$ for some $j\in\{k'+1,\ldots,n-k\}$,
then the all-zero choice of coefficients for rows
$\{k'+1,\ldots,k-1\}$ becomes illegal when giving row $k$ the
coefficient $1$, implying that
$|\calA|\leq\left(\frac{\delta}{2}+o(1)\right)2^k$, where $\delta =
1-2^{-(k-k')}$.
\item \label{item-rA+sA=0}If $M'_\calA$ is the $(k-1)\times n$ sub-matrix of rows $\{1,\ldots,k-1\}$ of
$M_\calA$ (that is, the matrix obtained by erasing the last row of
$M_\calA$), then every column of $M'_\calA$ contains at most 1
non-zero entry, and every row of $M'_\calA$ belongs to $\{0,\pm1\}^n
\setminus \{0,1\}^n$. To see this, notice that the coefficient of
row $k$ is set to $1$, otherwise we obtain at most $o(2^k)$ vectors.
We can thus regard this row as an affine vector in $\{0,1\}^n$, and
consider the $2^{k-1}$ combinations for the remaining rows. Now, a
column of $M'_\calA$ with at least $2$ non-zero entries implies that
the number of such legal combinations (resulting in a vector in
$\{0,1\}^n$) is at most $\frac{3}{4}\cdot 2^{k-1}$, and a row which
does not belong to $\{0,\pm1\}^n\setminus\{0,1\}^n$ implies that
this number is at most $2^{k-2}$. In both cases, we get $|\calA|
\leq (\frac{3}{8}+o(1))2^k$.
\item \label{item-M'-has-1-per-row}
Every row of $M'_\calA$ has at most $2$ non-zero values: assume that
the converse holds, that is, that row $m\in[k-1]$ contains at least
$2$ non-zero entries in indices $\{k+1,\ldots,n\}$. Since each of
the $k-1$ rows of $M'_\calA$ must contain a $-1$ value in an
exclusive column, it leaves at most $n-k-(k-1)=n-2k+1\leq 1$ column
(recall that $k \geq\frac{n}{2}$), which can contribute $1$
additional non-zero value to row $m$. We deduce that row $m$ has
precisely two non-zero entries at columns $\{k+1,\ldots,n\}$.
However, in this case column $m$ of $M_\calB$ has precisely two
non-zero entries , since \eqref{eq-n/2-MA-MB-structure-k+h=n} and
the orthogonality of $M_\calA,M_\calB$ imply that:
\begin{equation} \label{eq-MA-MB-orthogonality} (M_\calA)_{i,k+j} =
-(M_\calB)_{j,i} ~\mbox{ for all }i\in[k]\mbox{ and }j\in[h]
\end{equation} (the inner product of row $i$ of $M_\calA$ and row
$j$ of $M_\calB$ is $(M_\calA)_{i,k+j} + (M_\calB)_{j,i} = 0$). From
the same reason, column $k$ of $M_\calB$ has at least $k'$ non-zero
entries (as row $k$ of $M_\calA$ has the value $1$ in columns
$\{k+1,\ldots,k+k'\}$). Therefore, performing the process of Claim
\ref{clm-R-rows} first on column $m$ and then on column $k$ of
$M_\calB$ gives $|\calB| \leq \frac{3}{4}\cdot
\frac{2+o(1)}{\sqrt{\pi\ell}}$, hence the pair $\calA,\calB$ is
suboptimal.
\end{enumerate}
Items \eqref{item-rA+sA=0} and \eqref{item-M'-has-1-per-row} imply
that, if the pair $\calA,\calB$ is optimal, then without loss of
generality, $M'_\calA$ is of the form $\left(\begin{smallmatrix}
  I_{k-1} | 0 | -I_{k-1} | 0
\end{smallmatrix}\right)$, as each row has $1,-1$ in exclusive
columns and $0$ everywhere else. In particular, $k' = k-1$, and
since $k \geq n/2$ and $k+k' \leq n$, we get:
\begin{equation}\label{eq-rA=o(n)-opt1-k-bnd}
  k = h= \frac{n}{2}~~\mbox{ or }~~(k=\frac{n+1}{2}~,~h=\frac{n-1}{2})~,
\end{equation}
and without loss of generality (using the orthogonality of
$M_\calA,M_\calB$):
\begin{equation}\label{eq-rA=o(n)-opt1-MA-MB-struct}
   \mbox{\small$M_\calA = \left( \begin{array}{c|c|c||c}
  I_{k-1} & \begin{array}{c}
    0\\ \vdots \\ 0
  \end{array} & -I_{k-1} & \begin{array}{c}
    0\\ \vdots \\ 0
  \end{array} \\
  \hline
  0 & 1 & 1 \ldots 1 & 0/1
  \end{array}\right)~,~
  M_\calB = \left( \begin{array}{c|c|c||c}
  I_{k-1} & \begin{array}{c}
    -1\\ \vdots \\ -1
  \end{array} & I_{k-1} & \begin{array}{c}
    0\\ \vdots \\ 0
  \end{array} \\
  \hline\hline
  0 & 0/-1 & 0\ldots 0 & 1
  \end{array}\right)$} ~,\end{equation}
where the last column of $M_\calA$ and the last row and column of
$M_\calB$ do not exist in case $k=(n+1)/2$. If $h=n/2$ and
$(M_\calB)_{h,k}=0$ (as opposed to $-1$), then $|\calB| \leq
(1+o(1))2^h/\sqrt{\pi\ell}$: the first $h-1$ rows have at most
$(2+o(1))2^{h-1}/\sqrt{\pi\ell}$ combinations by the usual
Littlewood-Offord argument on column $k$, and when adding row $h$ we
must form an antichain. It follows that if $k=h=n/2$, then
$(M_\calB)_{h,k}=-1$ and, by orthogonality, $(M_\calA)_{k,n}=1$:
$$
   M_\calA = \left( \begin{array}{c|c|c||c}
  \multicolumn{3}{c||}{\ddots} & \begin{array}{c}\vdots \\ 0\end{array} \\
  \hline
  0 & 1 & 1 \ldots 1 & 1
  \end{array}\right)~,~
  M_\calB = \left( \begin{array}{c|c|c||c}
  \multicolumn{3}{c||}{\ddots}  & \begin{array}{c}\vdots \\ 0\end{array} \\
  \hline\hline
  0 & -1 & 0 \ldots 0& 1
  \end{array}\right)~.$$
Finally, notice that the above structure of $M_\calA$ implies that
the coefficient for row $k$ is always $1$: a coefficient of $0$
necessarily results in the all-zero vector, which is forbidden in
$\calA$ (for instance, since $|\calA|$ is an antichain, or since
$\ell > 0$). Therefore:
$$ |\calA| \leq 2^{k-1}~.$$
If $\chi_{B_1}^{(j)}\neq \chi_{B_1}^{(k+j)}$ for some $j \in [k-1]$,
we must assign the coefficient $0$ to row $j$ of $M_\calB$, and we
are done, as in this case $|\calB| \leq (1+o(1))2^h/\sqrt{\pi\ell}$.
Assume therefore that $\chi_{B_1}^{(j)}= \chi_{B_1}^{(k+j)}$ for all
$j \in [k-1]$, and define: $$P = \{i \in [h]: k+ i \notin B_1\}
=\{i\in[h] : \chi_{B_1}^{(k+i)}=0\}~,~Q=[h]\setminus P~.$$ Every row
$i\in P$ of $M_\calB$ has $\{0,1\}$ as the set of possible
coefficients, and every row $i \in Q$ has $\{0,-1\}$ as the possible
coefficients. Take $B\in\calB$, and suppose that the affine
combination which produces $B$ assigns the coefficient $1$ to $p$
rows of $P$ ($0\leq p\leq|P|$), and assigns the coefficient $-1$ to
$q$ rows of $Q$ ($0\leq q \leq |Q|$). It follows from
\eqref{eq-rA=o(n)-opt1-MA-MB-struct} that for all $A \in \calA$:
\begin{equation}\label{eq-ell-pq-relation}
\ell = |A \cap B| = p + (|Q|-q) + \chi_B^{(k)}~.\end{equation} Let
$\calB_0$ denote the sets $\{B\in \calB : k \notin B\}$, and let
$\calB_1 = \calB \setminus \calB_0$. By \eqref{eq-ell-pq-relation},
we obtain that $q = p + |Q| - \ell$ if $k \notin B$, hence:
$$ |\calB_0| \leq \sum_{p=0}^{|P|}
\binom{|P|}{p}\binom{|Q|}{p+|Q|-\ell} = \sum_{p=0}^{|P|}
\binom{|P|}{p}\binom{|Q|}{\ell - p} = \binom{h}{\ell}~.$$ Similarly,
if $k \in B$ then $q = p + |Q| - \ell + 1$, and it follows that:
$|\calB_1| \leq \binom{h}{\ell-1}$. Altogether:
$$|\calB| = |\calB_0|+|\calB_1| \leq  \binom{h}{\ell} + \binom{h}{\ell-1} = \binom{h+1}{\ell}~,$$
and as $|\calA| \leq 2^{k-1}$:
\begin{equation}\label{eq-h-ell-h-ell-1-bound} |\calA||\calB| \leq
\binom{h+1}{\ell} 2^{n-h-1}~.\end{equation} As the maxima of the
function $f(x)=\binom{x}{\ell}2^{-x}$ on the domain $\mathbb{N}$ are
achieved at $x\in\{2\ell-1,2\ell\}$, we conclude that $h \in
\{2\ell-2,2\ell-1\}$ (otherwise $|\calA||\calB| <
\binom{2\ell}{\ell}2^{n-2\ell}$). Finally, recalling that:
\begin{equation}
  \label{eq-chi=B-k-pq-relation}
  \chi_B^{(k)} = q - p + \chi_{B_1}^{(k)}~,
\end{equation}
and combining \eqref{eq-ell-pq-relation} and
\eqref{eq-chi=B-k-pq-relation} we get: $$\ell = |Q| +
\chi_{B_1}^{(k)} ~. $$ Therefore, whenever $\chi_{B_1}^{(k)}=0$ we
get $|Q|=\ell$, hence $B=  \cup_{i\in[\ell]} \{(i,k+i)\}$ for some
$B \in \calB$. Letting $B_1$ denote this set $B$ without loss of
generality, we obtain the statement of
\eqref{eq-rA-sA=o(n)-opt-family-1}.

Finally, let us link the above to the optimal family
\eqref{eq-opt-pair'}. Define: $$ X = \left\{\begin{array}
  {ll}\{k,n\} & \mbox{if }k=\frac{n}{2}\\
\{k\} & \mbox{if }k=\frac{n+1}{2}
\end{array}\right. ~.$$ Each set
$A\in\calA$ is obtained by choosing one out of each pair of elements
$\big\{\{i,k+i\}:i\in[k-1]\big\}$, then adding these $k-1$ chosen
elements to the elements of $X$. Define:
$$ Y = \left\{\begin{array}
  {ll}
\big\{\{i,k+i\}:i\in[k-1]\big\} \cup \big\{\{n\}\big\} & \mbox{if }k=\frac{n}{2}\\
\big\{\{i,k+i\}:i\in[k-1]\big\} & \mbox{if }k=\frac{n+1}{2}\\
\end{array}\right. ~.$$
Each set $B \in \calB_1$ (that is, those sets which contain $k$)
has, in addition to $k$, $\ell-1$ objects of $Y$. Each set $B \in
\calB_0$ is the union of $\ell$ objects of $Y$, and altogether, all
sets $B\in\calB$ are the union of $\ell$ objects of $Y \cup
\big\{\{k\}\big\}$. As the last set holds the $k-1$ pairs
$\{i,k+i\}$ for $i\in[k-1]$ and the single elements corresponding to
$X$, this fits the description of \eqref{eq-opt-pair'} for
$\kappa=h+1$, $\tau=k-1$ and swapping $\calA,\calB$.
\end{proof}

\subsection{The optimal family
\eqref{eq-rA-sA=o(n)-opt-family-2}}\label{sec::clm-O(1)-1-columns}
\begin{claim}\label{clm-O(1)-1-columns}
  If $O(1)$ of columns $\{k+1,\ldots,k+k'\}$ of $M_\calA$ are
  $\{1\}$-columns, then \eqref{eq-rA-sA=o(n)-opt-family-2} holds.
\end{claim}
\begin{proof} By the assumption and by
\eqref{eq-(1)-columns-and-(1,1)-columns}, we obtain that $k'-O(1)$
of the columns $\{k+1,\ldots,k+k'\}$ are $\{1,1\}$-columns, that is,
there are $k'-O(1)$ columns $j\in\{k+1,\ldots,k+k'\}$ where there
are precisely two non-zero entries in rows $\{k'+1,\ldots,k\}$, and
both entries are equal to $1$. For each such column $j$, let
$i_1(j),i_2(j) \in \{k'+1,\ldots,k\}$ denote the rows where these
two entries are located. Assume that, without loss of generality,
the pair of rows $k-1,k$ is the most popular pair among the above
pairs of rows $\{ (i_1(j),i_2(j)): j \mbox{ is a
$\{1,1\}$-column}\}$; it follows that there are $\omega(1)$ columns
(and in fact, $\Omega(k')$ columns) $j\in\{k+1,\ldots,k+k'\}$ such
that:
$$\left\{\begin{array}{l}
(M_\calA)_{k-1,j}=(M_\calA)_{k,j}=1~,\\
(M_\calA)_{i,j}=0 \mbox{ for all }i\in\{k'+1,\ldots,k-2\}~.
\end{array}\right.
$$
Hence, if we assign the same coefficient to rows $k-1,k$ then we
obtain $\omega(1)$ values which differ from $1$ in columns
$\{k+1,\ldots,k'\}$, and contribute $o(2^k)$ vectors to $\calA$. We
must therefore assign the coefficient $1$ to precisely one of the
rows $k-1,k$ (and assign the coefficient $0$ to the other).

The arguments given in the proof of Claim
\ref{clm-omega(1)-1-columns} regarding row $k$ readily imply the
following analogous results on rows $k-1,k$:
\begin{enumerate}[(1)]
\item  Rows $k-1,k$ contain the value $1$ in columns
$\{k+1,\ldots,k\}$.
\item Rows $k-1,k$ belong to $\{0,1\}^n$.
\item If $M'_\calA$ is the $(k-2)\times n$ sub-matrix of rows
$\{1,\ldots,k-2\}$ of $M_\calA$, then every column of $M'_\calA$
contains at most 1 non-zero entry, and every row of $M'_\calA$
belongs to $\{0,\pm1\}^n \setminus \{0,1\}^n$.
\item Every row of $M'_\calA$ contains at most $2$ non-zero entries.
\end{enumerate}
By the last two items, we deduce that if $\calA,\calB$ is an optimal
pair, then without loss of generality,
$M'_\calA=\left(\begin{smallmatrix}
  I_{k-2} | 0 | -I_{k-2} | 0
\end{smallmatrix}\right)$, and in particular, $k'=k-2$. The constraints
$k \geq n/2$ and $k+k'\leq n$ now imply:
\begin{equation}\label{eq-rA=o(n)-opt2-k-bnd}
  k = h= \frac{n}{2}~~\mbox{ or
  }~~(k=\frac{n+1}{2}~,~h=\frac{n-1}{2})
  ~~\mbox{ or }~~(k=\frac{n}{2}+1~,~h=\frac{n}{2}-1)~,
\end{equation}
and by orthogonality:
\begin{equation}
  \label{eq-rA=o(n)-opt2-MA-MB-struct}
 \mbox{\small$M_\calA = \left( \begin{array}{c|c|c|c||c|c}
  I_{k-2} & \begin{array}{c}
    0\\ \vdots \\ 0
  \end{array} & \begin{array}{c}
    0\\ \vdots \\ 0
  \end{array} & -I_{k-2} & \begin{array}{c}0 \\ \vdots \\ 0\end{array} & \begin{array}{c}0\\ \vdots \\ 0\end{array}\\
  \hline
  0 & 1 & 0 & 1 \ldots 1 & 0/1 & 0/1 \\
  \hline
  0 & 0 & 1 & 1 \ldots 1 & 0/1 & 0/1
  \end{array}\right) ,
    M_\calB = \left( \begin{array}{c|c|c|c||c|c}
  I_{k-2} & \begin{array}{r}
    -1\\ \vdots \\ -1
  \end{array} & \begin{array}{r}
    -1\\ \vdots \\ -1
  \end{array} & I_{k-2} & \begin{array}{c}0\\ \vdots \\ 0\end{array} &
\begin{array}{c}0\\ \vdots \\ 0\end{array}
   \\
  \hline\hline
  0 & 0/-1 & 0/-1 & 0\ldots 0 & 1 & 0\\
  \hline
    0 & 0/-1 & 0/-1 & 0\ldots 0 & 0 & 1\\
  \end{array}\right)$}~,
  \end{equation}
where the last two columns of $M_\calA$ and the last two rows and
columns of $M_\calB$ are optional, depending on whether
$k=\frac{n}{2}+1$, $k=\frac{n+1}{2}$ or $k=\frac{n}{2}$ (where we
have $0$, $1$ or $2$ of the last columns of $M_\calA$ and the last
rows and columns of $M_\calB$ respectively).

By \eqref{eq-rA=o(n)-opt2-MA-MB-struct}, it now follows that
choosing the same coefficient for both rows $k-1,k$ does not produce
sets in $\calA$ (so far we only showed that it produces $o(2^k)$
sets in $\calA$). Indeed, assigning the coefficient $0$ to both
these rows can only yield the all-zero vector, forbidden in $\calA$
(for instance, as $\ell > 0$). Assigning the coefficient $1$ to rows
$k-1,k$ can only yield a vector which is $1$ in every coordinate
$j\in[2k-2]$, and is the sum of the two rows $k-1,k$ in columns
$2k-1,2k$ if these columns exist. Hence, if this vector belongs to
$\{0,1\}^n$, then it contains any set which can be produced from
$M_\calA$, and we have $|\calA|=1$, and a suboptimal pair
$\calA,\calB$. It follows that: $$ |\calA| \leq 2^{k-1} ~.$$

Our next goal is to show that if row $q\in\{k-1,k\}$ of $M_\calB$
exists, then its entries in columns $k-1,k$ (marked by $0/-1$ in
\eqref{eq-rA=o(n)-opt2-MA-MB-struct}) are both $-1$. Let $q \in
\{k-1,k\}$ denote a row of $M_\calB$, let $m\in\{1,2\}$ denote the
number of rows of $\{k-1,k\}$ in $M_\calB$, and let $q' \neq q$
denote the additional row of $\{k-1,k\}$ in $M_\calB$ if $m=2$.
Since $m=1$ iff $k=\frac{n+1}{2}$ and $m=2$ iff $k=\frac{n}{2}$, it
follows that $m=2-(k-h)$.

First, assume that $\mathbf{(M_\calA)_{q,k-1}=(M_\calA)_{q,k}=0}$.
It follows that row $q$ is in $\{0,1\}^n$, and since $\calB$ is an
antichain, we get an additional factor of $\frac{1}{2}$ on $|\calB|$
(first apply the Littlewood-Offord Lemma on the remaining rows with
respect to column $k$, then consider the coefficient for row $q$).
It follows that $|\calB| \leq (1+o(1))\frac{2^h}{\sqrt{\pi\ell}}$,
and that $|\calA||\calB| \leq (\frac{1}{2}+o(1))2^n/\sqrt{\pi\ell}$.

Second, assume that $\mathbf{(M_\calA)_{q,k-1} \neq
(M_\calA)_{q,k}}$. Let $t_1$ denote the number of sets $B\in\calB$
produced from $M_\calB$ by assigning the coefficient $\alpha \neq 0$
to row $q$, and the coefficient $0$ to row $q'$ (if this row
exists), and let $t_2 = |\calB| - t_1$. Consider a set $B$ counted
by $t_2$: since row $q'$ does not take part in the affine
combinations, the combination of rows $[k-2]$ together with
$\chi_{B_1}$ sums up to the same value, some $\lambda$, in the two
columns $k-1,k$ (these two columns are identical in rows $[k-2]$).
The fact that indices $k-1,k$ of the resulting vector, $\chi_B$, are
$\{\lambda,\lambda-\alpha\}$, forces $\lambda$ to be equal to
$\alpha$. We can thus apply the Littlewood-Offord Lemma on rows
$[k-2]$ (with respect to column $k$, which has $1$ target value),
and deduce that:
$$t_1 \leq (1+o(1))\frac{2^{k-2}}{\sqrt{\pi\ell}}~.$$
To obtain an upper bound on $t_2$, for each of the remaining $2^m-1$
combinations of rows $\{k-1,k\}$ in $M_\calB$, column $k$ has at
most $2$ target values (in order to give a $\{0,1\}$ final value),
hence, by the Littlewood-Offord Lemma:
$$ t_2 \leq (2^m-1)(2+o(1))\frac{2^{k-2}}{\sqrt{\pi\ell}} ~.$$
It follows that:
$$ |\calB| = t_1 + t_2 \leq (2-2^{-m}+o(1))\frac{2^{m+k-2}}{\sqrt{\pi\ell}} =
(2 - 2^{-m}+o(1)) \frac{2^h}{\sqrt{\pi\ell}}~,$$ where in the last
equality we used the fact that $m=2-(k-h)$. The fact that $|\calA|
\leq 2^{k-1}$ now implies that the pair $\calA,\calB$ is suboptimal.

Having ruled out the cases $(M_\calA)_{q,k-1}=(M_\calA)_{q,k}=0$ and
$(M_\calA)_{q,k-1}\neq(M_\calA)_{q,k}$, we deduce that:
$$(M_\calA)_{q,k-1}=(M_\calA)_{q,k}=-1~,$$ hence the structure
of $M_\calA,M_\calB$ is:
$$\mbox{\small$M_\calA = \left( \begin{array}{c|c|c|c||c|c}
  \multicolumn{4}{c||}{\ddots} & \begin{array}{c}\vdots \\ 0\end{array} & \begin{array}{c}\vdots \\ 0\end{array}\\
  \hline
  0 & 1 & 0& 1 \ldots 1 & 1 & 1 \\
    \hline
  0 & 0 & 1& 1 \ldots 1 & 1 & 1
  \end{array}\right)~,~
  M_\calB = \left( \begin{array}{c|c|c|c||c|c}
  \multicolumn{4}{c||}{\ddots}  & \begin{array}{c}\vdots \\ 0\end{array}
   & \begin{array}{c}\vdots \\ 0\end{array} \\
  \hline\hline
  0 & -1 &-1 & 0 & 1 & 0 \\
  \hline
    0 & -1 &-1 & 0 & 0 & 1
  \end{array}\right)$}~,$$  as specified in
\eqref{eq-rA-sA=o(n)-opt-family-2}. To conclude the proof of the
claim, recall that every $A\in\calA$ has precisely one of the
elements $k-1,k$, hence the analysis of $|A\cap B|$ for all
$B\in\calB$ is exactly the same as in Claim
\ref{clm-omega(1)-1-columns} (precisely one of the columns $k-1,k$
of $M_\calB$ effects the intersection). It follows that
$|\calA||\calB| \leq
\left(\binom{h}{\ell}+\binom{h}{\ell-1}\right)2^{n-h-1}=\binom{h+1}{\ell}2^{n-h-1}$,
and hence $h \in \{2\ell - 2,2\ell-1\}$, otherwise $\calA,\calB$ is
a suboptimal pair. Similarly, the arguments of Claim
\ref{clm-omega(1)-1-columns} imply that $|Q|=\ell$, where $Q$ is the
set of indices $\{i\in[h]: k+i \in B_1\}$, and without loss of
generality, we can take $B_1$ to be $\cup_{i\in[\ell]}\{ i,k+i \}$.
Altogether, \eqref{eq-rA-sA=o(n)-opt-family-2} holds.

It remains to link the above to the optimal family
\eqref{eq-opt-pair'}. Define: $$ X = \left\{\begin{array}
  {ll}\{n-1,n\} & \mbox{if }k=\frac{n}{2}\\
  \{n\} & \mbox{if }k=\frac{n+1}{2}\\
\emptyset & \mbox{if }k=\frac{n}{2}+1
\end{array}\right. ~.$$ Recall that precisely one of the rows
$k-1,k$ receives the coefficient $1$ in a linear combination which
produces some $A\in\calA$ from $M_\calA$. It follows that each set
$A\in\calA$ is obtained by choosing one out of each pair of elements
$\big\{\{i,k+i\}:i\in[k-2]\big\} \cup \big\{ \{k-1,k\} \big\}$, then
adding these $k-1$ chosen elements to the elements of $X$. Define:
$$ Y = \left\{\begin{array}
  {ll}
\big\{\{i,k+i\}:i\in[k-2]\big\} \cup \big\{ \{n-1\},\{n\}\big\} & \mbox{if }k=\frac{n}{2}\\
\big\{\{i,k+i\}:i\in[k-2]\big\} \cup \big\{ \{n\}\big\} & \mbox{if }k=\frac{n+1}{2}\\
\big\{\{i,k+i\}:i\in[k-2]\big\} & \mbox{if }k=\frac{n}{2}+1\\
\end{array}\right. ~.$$
Recall that, for all $B \in\calB$, the elements $k-1,k$ are either
both in $B$ or both not in $B$. If $k-1,k \notin B$, then $B$ is the
union of $\ell$ elements of $Y$. Otherwise, $B$ contains, in
addition to $\{k-1,k\}$, the union of $\ell-1$ elements of $Y$.
Altogether, all sets $B\in\calB$ are the union of $\ell$ objects of
$Y \cup \big\{\{k-1,k\}\big\}$. As the last set holds the $k-2$
pairs $\{i,k+i\}$ for $i\in[k-2]$, the pair $\{k-1,k\}$ and the
single elements corresponding to $X$, this fits the description of
\eqref{eq-opt-pair'} for $\kappa=h+1$, $\tau=k-1$ and swapping
$\calA,\calB$.

This completes the proof of Claim \ref{clm-O(1)-1-columns} and of
Lemma \ref{lem-rA=o(n)}.
\end{proof}

\section{Proof of Lemma \ref{lem-rA+sA=Omega(n)}}\label{sec::lem-rA+sA=Omega(n)}
The assumption that $r_A + s_A = \Omega(n)$ implies that $|\calA| =
O(2^k / \sqrt{n})$. Thus, if $r_B + s_B = \omega(1)$ we deduce that
$|\calA||\calB| = o(2^n / \sqrt{n})$ and we are done. Assume
therefore that $r_B + s_B = O(1)$, and let $C_2 = [h] \setminus (R_B
\cup S_B)$. By definition of $R_B$ and $S_B$, the following holds:
\begin{itemize}
\item Every column of $M_\calB$ contains at most $1$ non-zero value
in the rows of $C_2$.
\item Every row of $C_2$ belongs to $\{0,\pm1\}^n \setminus
\{0,1\}^n$.
\end{itemize}
We wish to show that $M_\calB$ is roughly of the form $\left(-I_h
\mid 0 \mid I_h \right)$, although so far we did not obtain any
restriction on the number of rows in $C_2$ with more than $2$
non-zero entries in $M_\calB$. In contrast to the analysis of
$M_\calA$ in Lemma \ref{lem-rA=o(n)}, this does not follow directly
from the fact that $r_B+s_B=O(1)$, as $h$ might be substantially
smaller than $n/2$ (as opposed to $k$).

We therefore return to $M_\calA$ and claim that at most $O(1)$
columns of $M_\calA$ contain at least $2$ non-zero entries in a
\textbf{cascading} manner. In other words, the process where we
repeatedly select an arbitrary column of $M_\calA$ with at least two
non-zero entries, and remove the rows where it is non-zero from the
matrix, ends after at most $O(1)$ steps. To see this, assume that
$\omega(1)$ such columns exist: $j_1,\ldots,j_m$. Perform the
process of creating $R_A$, beginning with the above columns: choose
column $j_i$ at step $i$ for $i \leq m$, and complete the process in
an arbitrary order of column selection, $j_{m+1},\ldots,j_t$. By the
assumption of the lemma, $r_A + s_A = \Omega(n)$, hence two cases
are possible:
\begin{itemize} \item $r_A=o(n)$: in this case $s_A = \Omega(n)$.
Clearly, $r_A \geq 2m = \omega(1)$ by the assumption, and the
additional $O(1/\sqrt{n})$ factor resulting from the rows $S_A$
implies that $|\calA| = o(2^k / \sqrt{n})$.
\item $r_A=\Omega(n)$: by definition, $r_A = \sum_{i=1}^t r_i$. If
for some $i,j \leq t$ we have $r_i,r_j = \omega(\sqrt{n})$ then
$|\calA| = o(2^k / \sqrt{n})$. Recall that if $t \geq 2 \log n$,
then $|\calA| \leq 2^k / (\frac{3}{4})^t \leq 2^k / n$. These two
facts imply that precisely one $i$ satisfies $r_i = \Omega(n)$.
Therefore, column $i$ gives a factor of $O(1/\sqrt{n})$, and the
remaining $t-1$ columns give a factor of $o(1)$ as $t \geq m =
\omega(1)$ and each such column contributes a factor of at most
$\frac{3}{4}$. Altogether, we deduce that $|\calA| = o(2^k /
\sqrt{n})$.
\end{itemize}
\begin{figure}
\centering \fbox{\includegraphics{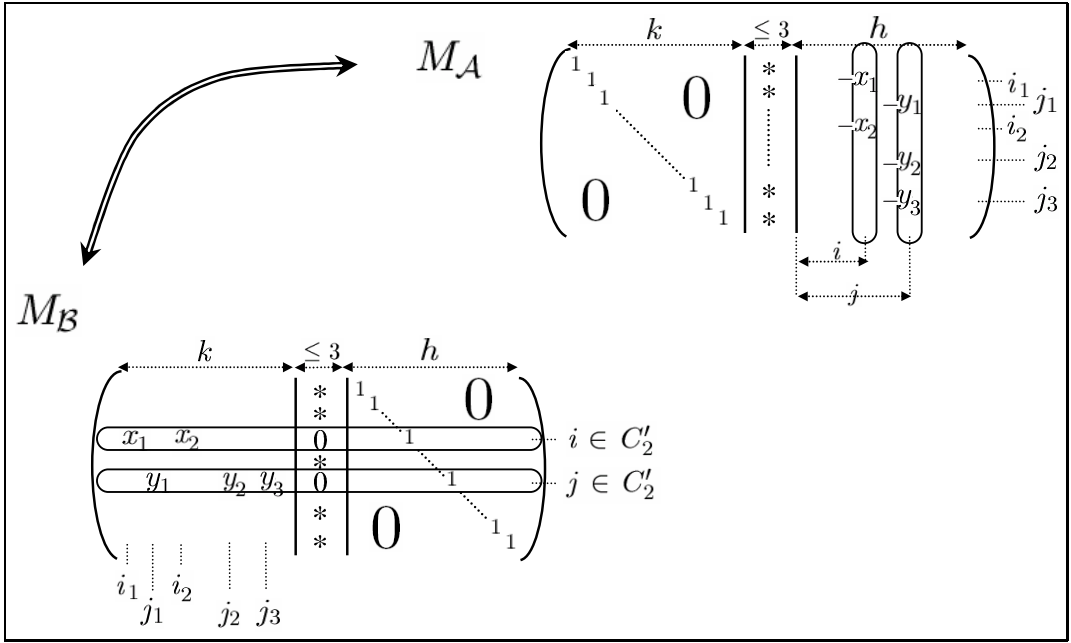}} \caption{The
duality between $M_\calA$ and $M_\calB$ when selected rows of
$M_\calB$ have 0 entries in columns $\{k+1,\ldots,n-h\}$.}
\label{fig::duality}
\end{figure}
Assume therefore that $M_\calA$ contains at most $O(1)$ columns
which contain at least $2$ non-zero entries in a cascading manner.
As we next show, returning to $M_\calB$, this implies that at most
$O(1)$ rows of $C_2$ contain more than $2$ non-zero entries. First,
recall that $k+h \geq n-3$ and that each column contains at most one
non-zero value in the rows of $C_2$. Thus, we can remove at most $3$
rows from $C_2$ and obtain a set $C_2'$, each remaining row of which
does not contain non-zero entries in indices $k+1,\ldots,n-h$.
Second, suppose rows $i,j\in C_2'$ each contains more than $2$
non-zero entries. Let $i_1,\ldots,i_r \in [k]$, $r \geq 2$, denote
the indices of the non-zero entries of row $i$ excluding its value
of $1$ at index $n-h+i$ (recall that columns $n-h+1,\ldots,n$ of
$M_\calB$ form the identity matrix of order $h$). Similarly, let
$j_1,\ldots,j_m \in [k]$, $m \geq 2$, denote the corresponding
indices of row $j$ :
$$ (M_\calB)_{i,i_t} \neq 0 \mbox { for }1\leq t \leq
r~,~(M_\calB)_{i,n-h+i} = 1~,$$
$$ (M_\calB)_{j,j_t} \neq 0 \mbox { for }1\leq t \leq
m~,~(M_\calB)_{j,n-h+j} = 1~.$$ Since the rows of $M_\calA$ are
orthogonal to the rows of $M_\calB$, and columns $1,\ldots,k$ of
$M_\calA$ form the identity matrix of order $k$, we deduce that:
$$ (M_\calA)_{n-h+i_t,i} \neq 0 \mbox { for }1\leq t \leq
r~,$$
$$ (M_\calA)_{n-h+j_t,j} \neq 0 \mbox { for }1\leq t \leq
m~.$$ See Figure \ref{fig::duality} for an illustration of the above
relation between $M_\calA$ and $M_\calB$. As the sets
$\{i_1,\ldots,i_r\}$ and $\{j_1,\ldots,j_m\}$ are disjoint, columns
$n-h+i$ and $n-h+j$ of $M_\calA$ each contains at least $2$ non-zero
entries in pairwise distinct indices. In general, if $m$ rows in
$C_2'$ contain more than $2$ non-zero entries, we deduce that $m$
columns in $M_\calA$ contain at least $2$ non-zero entries in a
cascading manner. As argued above, there are at most $O(1)$ such
columns in $M_\calA$, hence $m = O(1)$: let $C_2''$ denote the set
$C_2'$ after removing these $m$ rows, and let $h' = |C_2''| = h -
O(1)$. Each row of $C_2''$ is in $\{0,\pm1\}^n\setminus\{0,1\}^n$
and contains at most $2$ non-zero values, and we deduce that without
loss of generality:

\begin{equation}\label{eq-small-h-B-structure}
\begin{array}{c|c|c|c|c}
\multicolumn{2}{r@{\mbox{\tiny$ \dashrightarrow|$}}}{\mbox{\tiny$
\dashleftarrow \cdots \; k \; \cdots $}} &
\multicolumn{1}{@{\mbox{\tiny$|$}}c@{\mbox{\tiny$|$}}}{\mbox{\tiny$\dashleftarrow
\leq3\dashrightarrow$}} &
\multicolumn{2}{@{\mbox{\tiny$|\dashleftarrow$}}l}{\mbox{\tiny$\cdots \quad h \quad \cdots\dashrightarrow$}} \\
M_\calB = \left(\begin{array}{c}-I_{h'}\\ * \end{array}\right. &
\begin{array}{c}0\\ *
\end{array} &
\begin{array}{c}0\\ * \end{array} &
\begin{array}{c}I_{h'}\\
0
\end{array} &
\left.\begin{array}{c}0\\ I_{h-h'} \end{array}
\right)~.\end{array}\end{equation}
 Since the rows of $M_\calA$ and $M_\calB$ are
orthogonal, it follows that:
\begin{equation}\label{eq-small-h-A-structure}
\begin{array}{c|c|c|c|c}
\multicolumn{2}{r@{\mbox{\tiny$ \dashrightarrow|$}}}{\mbox{\tiny$
\dashleftarrow \cdots \quad k \quad \cdots $}} &
\multicolumn{1}{@{\mbox{\tiny$|$}}c@{\mbox{\tiny$|$}}}{\mbox{\tiny$\dashleftarrow
\leq3\dashrightarrow$}} &
\multicolumn{2}{@{\mbox{\tiny$|\dashleftarrow$}}l}{\mbox{\tiny$\cdots h \cdots\dashrightarrow$}} \\
M_\calA = \left(\begin{array}{c}I_{h'}\\ 0 \end{array}\right. &
\begin{array}{c}0\\ I_{k-h'}
\end{array} &
\begin{array}{c}*\\ * \end{array} &
\begin{array}{c}I_{h'}\\
0
\end{array} &
\left.\begin{array}{c}*\\ * \end{array}
\right)~.\end{array}\end{equation}
 The above structure of $M_\calA$
and $M_\calB$ provides an upper bound on $\ell$ in terms of $k$,
which we prove in Subsection \ref{sec::clm-ell-leq-k-2}:
\begin{claim}\label{clm-ell-leq-k/2}
Let $\calA$ and $\calB$ be as above. If $|\calA||\calB| =
\Omega(2^n/\sqrt{n})$, then: \begin{equation}\label{eq-ell-leq-k/2}
\ell \leq \left(\frac{1}{2}+o(1)\right) k~.
\end{equation}
\end{claim}
The proof of the lemma is completed by the next two claims, which
are proved in Subsections \ref{sec::clm-optimal-families} and
\ref{sec::clm-opt-rA=Omega(n)}:
\begin{claim}\label{clm-optimal-families}
  Let $\calA$ and $\calB$ be as above. If $r_A = o(n)$ then
  $|\calA||\calB| \leq \binom{2\ell}{\ell} 2^{n-2\ell}$.
  Furthermore, equality holds iff \eqref{eq-rA-sA=Omega(n)-opt-family}
  holds.
\end{claim}
\begin{claim}\label{clm-rA=Omega(n)}
  Let $\calA$ and $\calB$ be as above. If $r_A = \Omega(n)$ then
  the pair $\calA,\calB$ is suboptimal.
\end{claim}
\subsection{Proof of Claim
\ref{clm-ell-leq-k/2}}\label{sec::clm-ell-leq-k-2} Fix a choice of
coefficients for the rows $h'+1,\ldots,h$ of $M_\calB$, and let
$w_B$ denote the result of adding this combination to $\chi_{B_1}$.
As argued in the proof of Claim \ref{clm-n-geq-4ell}, the structure
of $M_\calB$ in \eqref{eq-small-h-B-structure} implies that each
index $j\in[h']$ such that
\begin{equation}
\label{eq-w-B-j-relation}w_B^{(j)}\neq 1- w_B^{(n-h+j)}
\end{equation} eliminates at least one of the two
possible coefficients for the row $j$ of $M_\calB$ (compare this to
the treatment of the vector $w_A$ in \eqref{eq-w-A-j-relation}).
Thus, if there are $\omega(1)$ such coefficients, then $w_B$ allows
at most $o(2^{h'})$ combinations of the remaining rows of $M_\calB$
to produce sets in $\calB$. Since $|\calA| = O(2^k /\sqrt{n})$
(recall that $r_A+s_A=\Omega(n)$), summing over at most $2^{h-h'}$
combinations for such vectors $w_B$ gives $o(2^{k+h}/\sqrt{n})$
pairs $(A,B)\in \calA \times \calB$.

It remains to treat vectors $w_B$ in which at most $O(1)$ indices
$j\in[h']$ satisfy \eqref{eq-w-B-j-relation}. Note that each $B \in
\calB$ produced from $w_B$ and a combination of rows $1,\ldots,h'$
of $M_\calB$ satisfies: \begin{equation}\label{eq-B-1-out-of-2-form}
|B \cap \{j,n-h+j\}| = 1 \mbox{ for all but at most }O(1)\mbox{
indices }j\in[h']~.\end{equation} Let $A\in\calA$, and let
$X_i\in\{0,1\}$ denotes the coefficient of the row $i$ of $M_\calA$
in the linear combination which produces $A$. By
\eqref{eq-B-1-out-of-2-form} and the structure of $M_\calA$ in
\eqref{eq-small-h-A-structure}, we obtain that: \begin{equation}
\label{eq-ell-part-1} \left|A \cap B \cap
\big([h']\cup\{n-h+1,\ldots,n-h+h'\}\big)\right| =
(\sum_{i=1}^{h'}X_i) + O(1)~.\end{equation}
 Furthermore, the structure of $M_\calA$ in
\eqref{eq-small-h-A-structure} gives:
\begin{equation}
  \label{eq-ell-part-2} \left|A \cap B \cap \{h'+1,\ldots,k\}\right| \leq \left|A \cap
\{h'+1,\ldots,k\}\right| = \sum_{i=h'+1}^k X_i~.
\end{equation}
Combining \eqref{eq-ell-part-1} and \eqref{eq-ell-part-2} with the
fact that $k+h'=n-O(1)$, we obtain that: $$\ell = |A\cap B|\leq
\bigg(\sum_{i=1}^k X_i \bigg) + O(1)~.
$$
Let $\epsilon > 0$, and assume that $\ell >
(1+\epsilon)\frac{k}{2}$. By the Chernoff bound, the number of
assignments of $\{0,1\}$ to the variables $X_1,\ldots,X_k$, which
satisfy $\sum_{i=1}^k X_i
> (1+\epsilon)\frac{k}{2}$, is at most
$2^k/\exp(\Omega(k))=2^k/\exp(\Omega(n))$. Therefore, the assumption
on $\ell$ implies that at most $O(2^k / \exp(\Omega(n))$ sets $A
\in\calA$ satisfy $|A\cap B|=\ell$, and summing over all sets $B$
whose vector $w_B$ is as above gives at most $2^{k+h} /
\exp(\Omega(n))$ pairs $(A,B)\in\calA\times\calB$. This contradicts
the assumption that $|\calA||\calB|=\Omega(2^n/\sqrt{n})$, and we
conclude that $\ell \leq (\frac{1}{2}+o(1))k$, as required. \qed

\subsection{Proof of Claim \ref{clm-optimal-families}}
\label{sec::clm-optimal-families} The assumptions $r_A + s_A =
\Omega(n)$ and $r_A = o(n)$ imply that $s_A = \Omega(n)$, and, as
before, we may assume that $r_A = O(1)$, otherwise we get $|\calA| =
o(2^k / \sqrt{n})$, leading to a suboptimal pair $\calA,\calB$.
Thus, each column of $M_\calA$ has at most $O(1)$ non-zero entries.
Since $n-(k+h) \leq 3$ and $h-h'=O(1)$, it follows that at most
$O(1)$ rows of $M_\calA$ have non-zero entries in columns
$\{k+1,\ldots,n-h\}\cup\{n-h+h'+1,\ldots,n\}$. Without loss of
generality, reorder the indices of these rows to $k'+1,\ldots,k$
(where $k' = k-O(1)$), and let $h'' = h' - O(1)$ reflect the
reordering of rows whose original indices belonged to $[h']$.
We obtain that:
\begin{equation}\label{eq-small-h-A-structure-rA=o(n)}
\begin{array}{c|c|c|c|c|c}
\multicolumn{3}{r@{\mbox{\tiny$ \dashrightarrow|$}}}{\mbox{\tiny$
\dashleftarrow \cdots \cdots \cdots \; \quad k \quad \; \cdots
\cdots \cdots $}} &
\multicolumn{1}{@{\mbox{\tiny$|$}}c@{\mbox{\tiny$|$}}}{\mbox{\tiny$\dashleftarrow
\leq 3\dashrightarrow$}} &
\multicolumn{2}{@{\mbox{\tiny$|\dashleftarrow$}}l}{\mbox{\tiny$\cdots h \cdots\dashrightarrow$}} \\
M_\calA =\left(\begin{array}{c}I_{h''}\\ 0 \\
0\end{array}\right. &
\begin{array}{c}0\\ I_{k'-h''}\\0
\end{array} &
\begin{array}{c}0\\ 0 \\ I_{k-k'}
\end{array} &
\begin{array}{c}0\\ 0\\ * \end{array} &
\begin{array}{c}I_{h''}\\0\\ 0
\end{array} &
\left.\begin{array}{c}0\\ 0\\ * \end{array}\right)~, \end{array}
\end{equation} and by the orthogonality of $M_\calA$ and
$M_\calB$:
\begin{equation}\label{eq-small-h-B-structure-rA=o(n)}
\begin{array}{c|c|c|c|c}
\multicolumn{2}{r@{\mbox{\tiny$\dashrightarrow|$}}}{\mbox{\tiny$
\dashleftarrow \cdots \; k' \; \cdots $}} &
\multicolumn{1}{@{\mbox{\tiny$|$}}c@{\mbox{\tiny$|$}}}{\mbox{\tiny$\dashleftarrow
O(1)\dashrightarrow$}} &
\multicolumn{2}{@{\mbox{\tiny$|\dashleftarrow$}}l}{\mbox{\tiny$\cdots\quad h \quad\cdots\dashrightarrow$}} \\
M_\calB =\left(\begin{array}{c}-I_{h''}\\ 0 \end{array}\right. &
\begin{array}{c}0\\ 0
\end{array} &
\begin{array}{c}0\\ * \end{array} &
\begin{array}{c}I_{h''}\\
0
\end{array} &
\left.\begin{array}{c}0\\ I_{h-h''} \end{array}
\right)~.\end{array}\end{equation} Notice that the first $k'$ rows
of $M_\calA$ form an antichain on the first $k'$ elements, hence:
$$|\calA| \leq (1+o(1))\frac{2^k}{\sqrt{\pi k' / 2}} \leq
(1+o(1))\frac{2^k}{\sqrt{\pi \ell}}~,$$ where the last inequality is
by \eqref{eq-ell-leq-k/2}. This yields an upper bound on
$|\calA||\calB|$ which is asymptotically tight, hence any additional
constant factor bounded away from $1$ which multiplies either
$|\calA|$ or $|\calB|$ implies that the pair $(\calA,\calB)$ is
suboptimal. In particular:
\begin{enumerate}[(i)]\item
\label{item-opt-1}If $k+h < n$, we have a suboptimal pair:
$|\calA||\calB| \leq
\left(\frac{1}{2}+o(1)\right)2^n/\sqrt{\pi\ell}$. Assume therefore
that $k+h=n$.
\item \label{item-opt-2}If $M_\calB$ has a column with more than $1$ non-zero entry,
we gain a multiplicative factor of at most $\frac{3}{4}$ and we are
done. The same applies to $M_\calA$: such a column has $O(1)$
non-zero entries, as $r_A = O(1)$, and once we set the combination
of these rows (gaining a factor of at most $\frac{3}{4}$) as well as
of rows $k'+1,\ldots,k$, the remaining $k'-O(1)$ rows out of $[k']$
must still form an antichain.
\item \label{item-opt-3}If $M_\calA$ has a row with more than $2$ non-zero entries, by Item \eqref{item-opt-1}
it corresponds to a column with more than $1$ non-zero entry in
$M_\calB$ (since statement \eqref{eq-MA-MB-orthogonality} holds),
which does not exist according to Item \eqref{item-opt-2}. The same
applies to the rows of $M_\calB$.
\item \label{item-opt-4} Each row of $M_\calB$ must belong
to $\{0,\pm1\}^n\setminus \{0,1\}^n$, otherwise the arguments of
Claim \ref{clm-S-rows} imply a constant multiplicative factor of at
most $\frac{1}{2}$.
\end{enumerate}
Items \eqref{item-opt-3} and \eqref{item-opt-4} imply that every row
of $M_\calB$ has precisely two non-zero entries: $\{1,-1\}$, and
without loss of generality, $h'' = h$. Recalling
\eqref{eq-small-h-A-structure-rA=o(n)} and
\eqref{eq-small-h-B-structure-rA=o(n)}, $M_\calA$ and $M_\calB$ take
the following form:
\begin{equation}\label{eq-sA=Omega(n)-MA-MB-form}
\begin{array}{lllll} M_\calA = \bigg( ~~\begin{array}{c}I_h\\
0
\end{array}
&\bigg|
\begin{array}{c}0\\ I_{k-h} \end{array}
&\bigg|
\begin{array}{c}I_{h}\\
0
\end{array}
&\bigg)~,~\\
\noalign{\medskip} M_\calB = \bigg( -I_h &\bigg| \quad\; 0 &\bigg|
\begin{array}{c}I_h\end{array} &\bigg)~.\end{array}
\end{equation}
Notice that the above structure of $M_\calB$ implies that
$\chi_B^{(j)}=\chi_{B_1}^{(j)}$ for all $j\in\{h+1,\ldots,k\}$ and
$B\in\calB$. As we assumed in \eqref{eq-A,B-contain-all-elements}
that $\bigcup_{B\in\calB}\calB=[n]$, it follows that
$\{h+1,\ldots,k\}\in B_1$.

Consider the rows of $M_\calB$, let $w_B$ take the initial value of
the vector $\chi_{B_1}$, then subtract from $w_B$ each row $i$ of
$M_\calB$ for which $k+i \in B_1$. This translates the possible
coefficients for each row $i$ of $M_\calB$ to $\{0,1\}$; hence, the
characteristic vector of every element of $\calB$ is a sum of $w_B$
with a sub-sum of the rows of $M_\calB$. First,
$w_B^{(j)}=\chi_{B_1}^{(j)}=1$ for all $j \in \{h+1,\ldots,k\}$.
Second, the structure of $M_\calB$ \eqref{eq-sA=Omega(n)-MA-MB-form}
implies that, if $w_B^{(j)} \neq 1$ for some $j\in[h]$, then row $j$
cannot be added to $w_B$ to yield a vector in $\{0,1\}^n$. Since
this leads to a suboptimal pair $(\calA,\calB)$ (of size at most
$(\frac{1}{2}+o(1))2^n/\sqrt{\pi\ell}$), we deduce that:
$$ w_B = \big( \overbrace{1\ldots 1}^k ~ \overbrace{0\ldots
0}^h \big) ~.$$ The structure of $M_\calB$
\eqref{eq-sA=Omega(n)-MA-MB-form} implies that for every
$B\in\calB$, $\chi_B$ is of the form: $$\chi_B=\big(
\overbrace{0/1\ldots 0/1}^h ~\overbrace{1\ldots 1}^{k-h}
~\overbrace{1/0\ldots 1/0}^h\big)~,$$ where precisely one index in
each of the pairs $\{(1,k+1),\ldots,(h,k+h)\}$ is equal to $1$ in
$\chi_B$. If $X_i \in \{0,1\}$ denotes the coefficient of row $i$ of
$M_\calA$ in a combination that produces some $A\in M_\calA$, it
follows from \eqref{eq-sA=Omega(n)-MA-MB-form} that $\ell = |A \cap
B| = \sum_{i=1}^k X_i$ for all $B \in \calB$. By the properties of
the binomial distribution, we deduce that $ |\calA| \leq
\binom{k}{\ell} $, and altogether:
$$ |\calA||\calB| \leq 2^{n-k}\binom{k}{\ell}~.$$ The expression above realizes the
bound \eqref{eq-final-cross-bound} iff either $k = 2\ell$ or $k =
2\ell -1$, hence the final structure of the optimal pair
$(\calA,\calB)$ is as described in Lemma \ref{lem-rA+sA=Omega(n)}.
\qed

\subsection{Proof of Claim
\ref{clm-rA=Omega(n)}}\label{sec::clm-opt-rA=Omega(n)} The
assumption $r_A = \Omega(n)$ implies that, unless $s_A=O(1)$, we get
$|\calA|=o(2^k/\sqrt{k})=o(2^k/\sqrt{n})$ as required. However, if
we remove the rows $R_A$ from $[k]$, \eqref{eq-small-h-A-structure}
implies that only the columns $\{k+1,\ldots,n-h\} \cup
\{n-h+h',\ldots,n\}$ can contribute $-1$ entries to the remaining
rows, and each column has at most $1$ non-zero entry in each of
these rows. Since $n-(k+h)\leq 3$ and $h-h'=O(1)$, we deduce that
$[k]-r_A - s_A = O(1)$, and altogether:
$$ r_A = k - O(1)~. $$
\begin{definition} A column of $M_\calA$ is called ``heavy'' if it
contains $k-O(1)$ non-zero entries.\end{definition} The next
argument shows that there exists a heavy column in $M_\calA$. There
are at most $O(1)$ columns which may contain more than $1$ non-zero
entry in $M_\calA$ (as columns $[k]$ and $\{n-h+1,\ldots,n-h+h'\}$
contain a single non-zero entry of $1$). Therefore, there exists
some column $q\in[n]$ of $M_\calA$ with $\Omega(r_A) = \Omega(k)$
non zero entries. If some other column has $\omega(1)$ non-zero
entries in a cascading manner, we obtain $|\calA| =
o(2^k/\sqrt{n})$, and we are done. We deduce the column $q$ has $r_A
- O(1) = k-O(1)$ non-zero entries, therefore column $q$ is heavy.
Applying the Littlewood-Offord Lemma to the $k-O(1)$ rows where
column $q$ is non-zero at, we obtain that:
\begin{equation}\label{eq-A-rA=Omega(n)-bound}|\calA| \leq (2+o(1))\frac{2^k}{\sqrt{\pi k/2}}
\leq (2+o(1)) \frac{2^k}{\sqrt{\pi\ell}}~,\end{equation} where the
last inequality is by \eqref{eq-ell-leq-k/2}.

Let $q$ denote a heavy column of $M_\calA$. Lemma
\ref{lem-2-part-antichain} enables us to eliminate the case where
all non-zero entries of $q$ are $\pm1$. To see this, assume the
converse, and let:
$$ U = \{ i \in [k] : (M_\calA)_{i,q}=1\}~,~ V = \{ i \in [k] :
(M_\calA)_{i,q}=-1\}~.$$ Recall that $|U|+|V|=k-O(1)$, and take
$\epsilon
> 0$. If $|U|\geq (\frac{1}{2}+\epsilon)k$, then Chernoff's bound
implies that the number of sub-sums of the rows $U \cup V$ which
give a value of $\{0,1\}$ in this column is at most
$2^k/\exp(\Omega(k))$. We deduce $|U|=(\frac{1}{2}+o(1))k$ and that
$|V|=(\frac{1}{2}+o(1))k$.

Set $m = n-(k+h) + (h-h') = O(1)$. For each possible set of values
$\underline{x}\in\{0,1\}^{m}$ for columns $\{k+1,\ldots,n-h\} \cup
\{n-h+h',\ldots,n\}$, the family of all sets $A\in\calA$ which
matches the pattern $\underline{x}$ in the above set of columns is
an antichain, and either $|A\cap V|=|A\cap U|$ or $|A\cap V| =
|A\cap U|-1$. Therefore, Lemma \ref{lem-2-part-antichain} implies
that $|\calA| = O(2^k / k)= O(2^k / n)$. We may therefore assume
that:
\begin{equation}
  \label{eq-q-not-pm1}
  \mbox{Every heavy column $q$ of $M_\calA$ satisfies
  $(M_\calA)_{i,q}\notin\{0,\pm1\}$ for some }i\in[k]~.
\end{equation}
This provides an upper bound on $|\calB|$:
\begin{equation}
  \label{eq-b-1/2-factor}
  |\calB| \leq 2^{n-k-1}~.
\end{equation}
The above bound follows immediately if $h<n-k$, so consider the case
$k+h=n$, and let $q$ denote a heavy column of $M_\calA$. By the
orthogonality of $M_\calA,M_\calB$, \eqref{eq-MA-MB-orthogonality}
holds, and \eqref{eq-q-not-pm1} now implies that $(M_\calB)_{q-k,i}
\notin \{0,\pm 1\}$ for some $i\in[k]$. In particular, row $q-k$ of
$M_\calB$ does not belong to $\{0,\pm1\}^n$, and hence $ |\calB|
\leq 2^{h-1}$ (as enumerating on the coefficients for rows
$[h]\setminus\{q-k\}$ of $M_\calB$ leaves at most one legal
coefficient for row $q-k$).

Combining \eqref{eq-b-1/2-factor} with
\eqref{eq-A-rA=Omega(n)-bound} yields an asymptotically tight upper
bound on $|\calA||\calB|$:
$$|\calA||\calB| \leq (1+o(1))\frac{2^n}{\sqrt{\pi k / 2}} \leq
(1+o(1))\frac{2^n}{\sqrt{\pi \ell}}~.$$ Let $\epsilon > 0$; if $k
\geq (2+\epsilon)\ell$, then the first inequality of the bound above
implies that the pair $\calA,\calB$ is suboptimal. Therefore, adding
this to \eqref{eq-ell-leq-k/2}, we may assume that:
\begin{equation}\label{eq-k-eqq-2-ell}
k = (2 + o(1))\ell~.
\end{equation}
Next, we wish to eliminate the case where some column $q$ has
$k-O(1)$ non-zero entries, all of which have the same sign. In this
case, let $Q = \{i : (M_\calA)_{i,q} \neq 0\}$. As all the entries
in rows $Q$ and column $q$ of $M_\calA$ have the same sign, only the
all-zero linear combination of these rows can produce the value $0$
at index $q$. Applying the Littlewood-Offord Lemma to the rows $Q$,
we obtain an upper bound on the number of combinations which produce
the value $1$, and altogether:
$$|\calA| \leq 2^{k-|Q|} (\binom{|Q|}{\lfloor|Q|/2\rfloor}+1) =
(1+o(1))\frac{2^k}{\sqrt{\pi\ell}}~,$$ where in the last inequality
we used the fact that $|Q|\geq (2+o(1))\ell$, as $|Q|=k-O(1)$. By
\eqref{eq-b-1/2-factor}, this implies that $|\calA||\calB| \leq
(\frac{1}{2}+o(1))2^n/\sqrt{\pi\ell}$, implying the statement of the
claim. We thus assume that: \begin{equation}   \label{eq-q-2-vals}
\mbox{Every heavy column $q$ of $M_\calA$ contains both positive and
negative entries}~.
\end{equation}
Using the last statement, we prove the next claim:
\begin{claim}\label{clm-a-l-d-lambda} Let $\lambda \in \{0,1\}$, $L\subset [k]$ and $d
> 0$, and let $q$ denote a heavy column of $M_\calA$. Define:
\begin{equation}
  \label{eq-A-d-lambda-def}
  \calA_{L,d,\lambda}^{(q)} = \{ A \in \calA: |A\cap L|=d~,
\chi_{A}^{(q)}=\lambda\}~.
\end{equation} If $d= (1+o(1))\ell$ and
$|L|\geq(1+o(1))\ell$ then:
\begin{equation} \label{eq-A-lambda-d-bound} |\calA_{L,d,\lambda}^{(q)}|
\leq \left(\frac{3}{4}+o(1)\right)\frac{2^k}{\sqrt{\pi\ell}}~.
\end{equation}
\end{claim}
\begin{proof}
Let $Q$ denote the indices of the rows in which column $q$ of
$|\calA|$ has a non-zero entry. Observe that if $Q \nsubseteq L$,
then the rows of $L$ have at most $\binom{|L|}{d}$ legal
combinations, and the remaining rows $[k]\setminus L$ have at most
$2^{k-|L|-1}$ legal combinations, as these rows contain non-zero
entries in column $q$, which must combine to a final value of
$\lambda$. Hence, in this case: $$ |\calA_{L,d,\lambda}^{(q)}| \leq
\frac{1}{2} \cdot 2^{k-|L|}\binom{|L|}{d} \leq \frac{1}{2}
2^{k-|L|}\binom{|L|}{\lfloor|L|/2\rfloor} =
\frac{1+o(1)}{2}\cdot\frac{2^k}{\sqrt{\pi|L|/2}} \leq
\left(\frac{1}{\sqrt{2}}+o(1)\right)\frac{2^k}{\sqrt{\pi\ell}}~,
$$
where the last inequality is by the fact that $|L|\geq(1+o(1))\ell$.
 Assume therefore that $Q \subset L$, and notice that,
as $|Q| = k-O(1)$ and $L \subset [k]$, then $|L|=k-O(1)$, and by
\eqref{eq-k-eqq-2-ell}:
$$|L|= (2+o(1))\ell = (2+o(1))d~.$$
Fix an enumeration on the coefficients of the rows $[k]\setminus L$,
and let $\mathcal{S} \subset 2^{L}$ denote the $d$-element subsets
of the rows of $L$ which extend this enumeration to elements of
$\calA_{L,d,\lambda}^{(q)}$. Let $j_1,j_2\in L$ be two indices such
that $(M_\calA)_{j_1,q}\neq(M_\calA)_{j_2,q}$ (such indices exist by
\eqref{eq-q-2-vals} and since $Q\subset L$), and define:
$$\mathcal{S}_0 = \left\{ S \subset [L]: |S|=d~,~|S\cap\{j_1,j_2\}|=1\right\}~.$$
Notice that, as $j_1\neq j_2$, the function
$f:\mathcal{S}_0\to\mathcal{S}_0$ which swaps $j_1,j_2$ is a
bijection, which satisfies the following property for all $S \in
\mathcal{S}_0$: at most one of the subsets $\{S,f(S)\}$ can belong
to $\mathcal{S}$. Furthermore, if $S$ is a random $d$-element set of
$L$, then:
$$ \Pr[S \in \mathcal{S}_0] =
\frac{2\binom{|L|-2}{d-1}}{\binom{|L|}{d}} = \frac{2d
(|L|-d)}{|L|(|L|-1)} = \frac{1}{2}+o(1)~,$$ and thus
$|\mathcal{S}_0| = (\frac{1}{2}+o(1))\binom{|L|}{d}$, and we deduce
that:
$$ |\mathcal{S}| \leq \binom{|L|}{d} - \frac{|\mathcal{S}_0|}{2} =
\left(\frac{3}{4}+o(1)\right)\binom{|L|}{d} ~.$$ Therefore:
$$ |\calA_{L,d,\lambda}^{(q)}| \leq
2^{k-|L|}|\mathcal{S}| \leq
\left(\frac{3}{4}+o(1)\right)\frac{2^k}{\sqrt{\pi|L|/2}} =
\left(\frac{3}{4}+o(1)\right)\frac{2^k}{\sqrt{\pi\ell}}~,
$$
as required.
\end{proof}
In order to deduce the claim from \eqref{eq-A-lambda-d-bound}, we
treat the two cases $k+h<n$ and $k+h=n$ in Claims
\ref{clm-subopt-k+h<n} and \ref{clm-subopt-k+h=n} below.
\begin{claim}\label{clm-subopt-k+h<n} Let $\calA,\calB$ be as above. If $k+h < n$, then
the pair $\calA,\calB$ is suboptimal.
\end{claim} \begin{proof}In this case, we may
assume that $k + h = n - 1$, otherwise
\eqref{eq-A-rA=Omega(n)-bound} implies that $|\calA||\calB| \leq
(\frac{1}{2}+o(1))2^n/\sqrt{\pi\ell}$. Recalling
\eqref{eq-small-h-B-structure} and \eqref{eq-small-h-A-structure},
we have: \begin{equation}\label{eq-MA-MB-1-missing-column}
\begin{array}{c|c|c|c|c}
\multicolumn{2}{r@{\mbox{\tiny$\dashrightarrow|$}}}{\mbox{\tiny$\dashleftarrow
\cdots\cdots k \cdots\cdots $}} &
\multicolumn{1}{@{\mbox{\tiny$|$}}c@{\mbox{\tiny$|$}}}{\mbox{\tiny$\dashleftarrow1\dashrightarrow$}}
&
\multicolumn{2}{@{\mbox{\tiny$|\dashleftarrow$}}l}{\mbox{\tiny$\cdots\; h \;\cdots \dashrightarrow$}} \\
 M_\calA = \left( \begin{array}{c}I_{h'}\\ 0
\end{array}\right. &
\begin{array}{c}0 \\ I_{k-h'} \end{array} &
\begin{array}{c}* \\ * \end{array} &
\begin{array}{c}I_{h'}\\ 0 \end{array} &
\left.\begin{array}{c} *\\ * \end{array}\right) \\
\noalign{\medskip}
\multicolumn{2}{r@{\mbox{\tiny$\dashrightarrow|$}}}{\mbox{\tiny$\dashleftarrow
\cdots\cdots k \cdots\cdots $}} &
\multicolumn{1}{@{\mbox{\tiny$|$}}c@{\mbox{\tiny$|$}}}{\mbox{\tiny$\dashleftarrow1\dashrightarrow$}}
&
\multicolumn{2}{@{\mbox{\tiny$|\dashleftarrow$}}l}{\mbox{\tiny$\cdots\;
h \;\cdots \dashrightarrow$}} \\
 M_\calB = \left(
\begin{array}{c}-I_{h'}\\ * \end{array} \right.&
\begin{array}{c}0 \\ * \end{array} &
\begin{array}{c}0 \\ * \end{array} &
\begin{array}{c}I_{h'}\\ 0 \end{array} &
\left.\begin{array}{c}0\\ I_{h-h'} \end{array} \right)
\end{array}~.
\end{equation}
Let $m=h-h'=O(1)$, and consider a choice of coefficients for rows
$h'+1,\ldots,h$ of $M_\calB$, yielding (together with $\chi_{B_1}$)
a vector $w_B$. First, by \eqref{eq-A-rA=Omega(n)-bound}, each of
the $2^{m}-1$ choices of coefficients such that
$(w_B^{(n-m+1)}\ldots w_B^{(n)}) \neq 0$ can each be completed to a
pair $(A,B)\in \calA \times \calB$, in at most $$2^{h-m} \cdot
(2+o(1))\frac{2^k}{\sqrt{\pi\ell}} =
(1+o(1))\frac{2^{n-m}}{\sqrt{\pi\ell}}$$ ways. Let $\calB_0$ denote
the sets $B \in \calB$ which can be produced from the remaining
combination for $w_B$ (the one for which
$w_B^{(n-m+1)}=\ldots=w_B^{(n)}=0$). In order to show that
$\calA,\calB$ is suboptimal, it is enough to show that:
\begin{equation}   \label{eq-B0-suff-condition}
|\calA||\calB_0| \leq (\alpha+o(1))
\frac{2^{n-m}}{\sqrt{\pi\ell}}\mbox{ for some }\alpha <
1~,\end{equation} since this would imply:
\begin{equation}\label{eq-B0-suff-condition-pf}|\calA||\calB| \leq
(2^m-1)(1+o(1))\frac{2^{n-m}}{\sqrt{\pi\ell}} + (\alpha+o(1))
\frac{2^{n-m}}{\sqrt{\pi\ell}} =
\left(1-\frac{1-\alpha}{2^m}+o(1)\right)\frac{2^n}{\sqrt{\pi\ell}}~.\end{equation}
If for some index $j \in [h']$ we have $w_B^{(n-h+j)} \neq
1-w_B^{(j)}$, then row $j$ of $M_\calB$ has at most one legal
coefficient, hence $|\calB_0|\leq 2^{h-m-1}$, and the same holds in
case $w_B \notin \{0,1\}^n$ (if $j\in\{h'+1,\ldots,n-h\}$ is such
that $w_B^{(j)}\notin \{0,1\}$, then $\calB_0=\emptyset$). As
$|\calA| \leq (2+o(1))\frac{2^k}{\sqrt{\pi\ell}}$ and $k+h< n$, it
follows that in the above two cases $|\calA||\calB_0|\leq
\left(\frac{1}{2}+o(1)\right)\frac{2^{n-m}}{\sqrt{\pi\ell}}~,$
satisfying \eqref{eq-B0-suff-condition} for $\alpha=\frac{1}{2}$.

Assume therefore that $w_B^{(n-h+j)}=1-w_B^{(j)}$ for all
$j\in[h']$, and that $w_B \in \{0,1\}^n$, and define: $$L=[h']\cup
\{h'+1 \leq i \leq k:w_B^{(i)}=1\}~.$$ Recalling that
$w_B^{(n-h+h'+1)}=\ldots=w_B^{(n)}=0$,
\eqref{eq-MA-MB-1-missing-column} implies that every $B$ produced
from $w_B$ satisfies:
\begin{equation}\label{eq-ell-sum-miss-dim} \ell = |A \cap B| =
\mathbf{1}_{\{k+1\in A\cap B\}} +\sum_{i \in L} X_i~,\end{equation}
for all $A\in\calA$, where $X_i \in \{0,1\}$ denotes the coefficient
for row $i$ in a combination which produces $A$ from $M_\calA$. We
may assume that $\calB_0\neq\emptyset$ (otherwise
\eqref{eq-B0-suff-condition} immediately holds), and by
\eqref{eq-ell-sum-miss-dim} we obtain that $|L|\geq \ell-1$, and in
particular, $|L|\geq (1+o(1))\ell$.

If column $k+1$ of $M_\calA$ has $o(k)=o(|L|)$ non-zero entries in
some rows $U$, fix an enumeration on the coefficients of these rows,
and let $L' = L \setminus U$, noting that $|L'|=(1-o(1))|L|\geq
(1-o(1))\ell$. The enumeration on the coefficients for the rows $U$
determines whether or not $k+1\in A\cap B$, and by
\eqref{eq-ell-sum-miss-dim}, this determines the value of
$\sum_{i\in L'}X_i$. Therefore, by the properties of the binomial
distribution, there are at most $\binom{|L'|}{\lfloor|L'|/2\rfloor}
\leq 2^{|L'|} / \sqrt{\pi |L'|/2}$ combinations for the coefficients
of the rows $L'$. We conclude that:
\begin{itemize}
  \item In case $|L| \geq (1-o(1))k$, recalling
\eqref{eq-k-eqq-2-ell}, we get $|\calA| \leq
(1+o(1))\frac{2^k}{\sqrt{\pi\ell}}$.
\item Otherwise, $k-|L|=\Omega(k)$, and after choosing a combination for
the rows $L'$, we are left with rows $[k]\setminus(L \cup U)$ which
contain $\Omega(k)$ non-zero entries in some heavy column $q$ of
$M_\calA$ (recall that each heavy column has $k-O(1)$ non-zero
entries). The Littlewood-Offord Lemma gives a factor of
$O(1/\sqrt{k})$ on the number of combinations for the remaining
rows, which, when multiplied by the previous factor of
$O(1/\sqrt{|L'|})=O(1/\sqrt{k})$ gives $|\calA| \leq
O(2^k/k)=O(2^k/\ell)$. In particular, we have $|\calA| \leq
(1+o(1))\frac{2^k}{\sqrt{\pi\ell}}$ (with room to spare).
\end{itemize}
Altogether, as $|\calB_0| \leq 2^{h-m} \leq 2^{n-m-k-1}$, in both
cases we obtain that \eqref{eq-B0-suff-condition} holds for $\alpha
= \frac{1}{2}$.

It remains to treat the case where column $k+1$ of $M_\calA$ has
$\Omega(k)$ non-zero entries; by the arguments in the beginning of
the proof of Claim \ref{clm-rA=Omega(n)}, it follows that column
$k+1$ is heavy. Therefore, recalling that $\calB_0\neq\emptyset$ and
using the definition \eqref{eq-A-d-lambda-def}, it follows that:
$$|\calA| = \left\{\begin{array}
  {cl}
  |\calA_{L,\ell,0}^{(k+1)}|+|\calA_{L,\ell,1}^{(k+1)}| &
   \mbox{if  }w_B^{(k+1)}=0\\
   \\
  |\calA_{L,\ell,0}^{(k+1)}|+|\calA_{L,\ell-1,1}^{(k+1)}| &
     \mbox{if  }w_B^{(k+1)}=1
\end{array}\right.~.$$
Applying Claim \ref{clm-a-l-d-lambda} (recall that $|L|\geq \ell -
1$) gives:
$$ |\calA| \leq 2 \cdot \left(\frac{3}{4}+o(1)\right)\frac{2^k}{\sqrt{\pi\ell}} =
\left(\frac{3}{2}+o(1)\right)\frac{2^k}{\sqrt{\pi\ell}}~,$$ and as
$|\calB_0| \leq 2^{h-m} \leq 2^{n-m-k-1}$,
\eqref{eq-B0-suff-condition} holds for $\alpha=\frac{3}{4}$, as
required.
\end{proof}

\begin{claim}\label{clm-subopt-k+h=n}Let $\calA,\calB$ be as above. If $k+h=n$, then
the pair $\calA,\calB$ is suboptimal. \end{claim}
\begin{proof} The proof will follow from arguments similar to those
in the proof of Claim \ref{clm-subopt-k+h<n}; the factor of
$\frac{1}{2}$ which followed from the case $k+h<n$ is replaced by
the duality between $M_\calA,M_\calB$ \eqref{eq-MA-MB-orthogonality}
when $k+h=n$. The assumption $k+h=n$ gives
\eqref{eq-small-h-B-structure} and \eqref{eq-small-h-A-structure}
the following form:
$$
\begin{array}{r|c|c|l}
\multicolumn{2}{r@{\mbox{\tiny$\dashrightarrow|$}}}{\mbox{\tiny$\dashleftarrow
\cdots\cdots k \cdots\cdots $}} &
\multicolumn{2}{@{\mbox{\tiny$|\dashleftarrow$}}l}{\mbox{\tiny$\cdots\; h \;\cdots \dashrightarrow$}} \\
 M_\calA = \left( \begin{array}{c}I_{h'}\\ 0
\end{array} \right.&
\begin{array}{c}0 \\ I_{k-h'} \end{array} &
\begin{array}{c}I_{h'}\\ 0 \end{array} &
\left.\begin{array}{c} *\\ * \end{array} \right) \\
\noalign{\medskip}
\multicolumn{2}{r@{\mbox{\tiny$\dashrightarrow|$}}}{\mbox{\tiny$\dashleftarrow
\cdots\cdots k \cdots\cdots $}} &
\multicolumn{2}{@{\mbox{\tiny$|\dashleftarrow$}}l}{\mbox{\tiny$\cdots\cdots
h \cdots\cdots \dashrightarrow$}} \\ M_\calB = \left(
\begin{array}{c}-I_{h'}\\ * \end{array} \right.&
\begin{array}{c}0 \\ * \end{array} &
\begin{array}{c}I_{h'}\\ 0 \end{array} &
\left.\begin{array}{c}0\\ I_{h-h'} \end{array} \right)
\end{array}~.
$$
Let $q\in[n]$ denote a heavy column of $M_\calA$; by the above
structure of $M_\calA$, we can assume without loss of generality
that $q = n$. Let $p\in[k]$ be such that $(M_\calA)_{p,n} \notin
\{0,\pm 1\}$ (such a $p$ exists by \eqref{eq-q-not-pm1}). Recall
that, as $k+h=n$, the orthogonality of $M_\calA,M_\calB$ implies
that \eqref{eq-MA-MB-orthogonality} holds, and thus $(M_\calB)_{h,p}
= -(M_\calB)_{p,n} \notin \{0,\pm 1\}$.

Consider the following set of rows of $M_\calB$:
$$ W = \left\{\begin{array}
  {ll} \{p\}\cup
\{h'+1,\ldots,h-1\} & \mbox{if }p\in[h']~,\\
 \{h'+1,\ldots,h-1\} & \mbox{otherwise}~.
\end{array}\right.
$$
Let $m=|W|$, and consider one of the $2^m-1$ choices of coefficients
for the rows $W$ of $M_\calB$, such that the sum of $\chi_{B_1}$ and
the resulting combination of these rows, satisfies $w_B^{(k+j)}\neq
0$ for some $j \in W$. Observe that $w_B$ allows at most one
coefficient for row $h$ of $M_\calB$, since all the remaining rows
$[h-1]\setminus W$ have $0$ entries at column $p$, whereas
$(M_\calB)_{h,p}\notin\{0,\pm1\}$. Therefore, by
\eqref{eq-A-rA=Omega(n)-bound}, each of the $2^m-1$ possibilities
for such vectors $w_B$ can produce at most:
$$ 2^{h-m-1} \cdot (2+o(1))\frac{2^k}{\sqrt{\pi\ell}} =
(1+o(1))\frac{2^{n-m}}{\sqrt{\pi\ell}} $$ pairs
$(A,B)\in\calA\times\calB$. Consider the remaining combination of
the rows $W$, satisfying $w_B^{(k+j)}=0$ for all $j \in W$, and let
$\calB_0$ denote the sets $B\in\calB$ which can be produced from
$w_B$. Using this notation, it is enough to show that
\eqref{eq-B0-suff-condition} holds, and the claim will follow from
the resulting calculation \eqref{eq-B0-suff-condition-pf}.

As before, the fact that $(M_\calB)_{h,p}\notin\{0,\pm1\}$ and that
the remaining rows $[h-1]\setminus W$ have $0$ entries in column
$p$, implies that there is at most one coefficient possible for row
$h$. If no coefficient for row $h$ is legal, we get
$\calB_0=\emptyset$ and \eqref{eq-B0-suff-condition} holds,
otherwise let $\tilde{w}_B$ denote the sum of $w_B$ with the
appropriate multiple of row $h$ of $M_\calB$. We are left with
$h-m-1$ rows of $M_\calB$ whose coefficients were not yet
determined: rows $[h-1]\setminus W = [h']\setminus\{p\}$.

If $\tilde{w}_B^{(j)} \neq 1-\tilde{w}_B^{(k+j)}$ for some
$j\in[h']\setminus \{p\}$ or $\tilde{w}_B \neq \{0,1\}^n$, we obtain
an additional factor of at most $\frac{1}{2}$ from one of the
remaining rows of $M_\calB$, and $|\calB_0| \leq 2^{h-m-2}$.
Combining this with \eqref{eq-A-rA=Omega(n)-bound} implies that
\eqref{eq-B0-suff-condition} holds for $\alpha = \frac{1}{2}$.
Assume therefore that $\tilde{w}_B^{(j)}=1-\tilde{w}_B^{(k+j)}$ for
all $j\in[h']\setminus\{p\}$ and that $\tilde{w}_B \in \{0,1\}^n$,
and define: $$L=[h']\setminus\{p\} \cup \left\{i \in
\{h'+1,\ldots,k\}\cup \{p\}:\tilde{w}_B^{(i)}=1\right\}~.$$ Since
every set $B$ produced from $\tilde{w}_B$ satisfies
 $|B \cap \{j,k+j\}|=1$ for all $j\in[h']\setminus\{p\}$ and
 $k+j\notin B$ for all $j\in W$, we deduce that, if $p\notin[h']$
 (in which case $W=\{h'+1,\ldots,h-1\}$):
\begin{equation}\label{eq-ell-sum-full-dim}\ell = |A \cap B| =
\mathbf{1}_{\{n\in A\cap B\}} +\sum_{i \in L} X_i~,\end{equation}
for all $A\in\calA$, where $X_i \in \{0,1\}$ denotes the coefficient
for row $i$ in a combination which produces $A$ from $M_\calA$. On
the other hand, if $p\in[h']$, then $p \in W$ and it follows that
$\tilde{w}_B^{(k+p)}=0$, and:
\begin{itemize}
\item If $\tilde{w}_B^{(p)}=0$, then $p\notin L$, and indeed, $X_p$ does not
contribute to $|A\cap B|$ for all $A\in \calA$ and $B$ produced by
$\tilde{w}_B$, as neither $p$ nor $k+p$ belong to $B$.
\item If $\tilde{w}_B^{(p)}=1$, then $p\in L$, and indeed $X_p$ contributes $1$ to
$|A\cap B|$ for all $A\in \calA$ and $B$ produced by $\tilde{w}_B$,
as $p\in B$ and $k+p\notin B$.
\end{itemize} We deduce that
\eqref{eq-ell-sum-full-dim} holds for $p\in[h']$ as-well. Recalling
that $\calB_0\neq\emptyset$ (otherwise \eqref{eq-B0-suff-condition}
immediately holds) \eqref{eq-ell-sum-full-dim} gives $|L|\geq
\ell-1$, and in particular, $|L|\geq (1+o(1))\ell$. Using the
definition \eqref{eq-A-d-lambda-def}, it follows that:
$$|\calA| = \left\{\begin{array}
  {cl}
  |\calA_{L,\ell,0}^{(n)}|+|\calA_{L,\ell,1}^{(n)}| &
   \mbox{if  }\tilde{w}_B^{(n)}=0\\
   \\
  |\calA_{L,\ell,0}^{(n)}|+|\calA_{L,\ell-1,1}^{(n)}| &
     \mbox{if  }\tilde{w}_B^{(n)}=1
\end{array}\right.~.$$
Applying Claim \ref{clm-a-l-d-lambda} (recall that $|L|\geq \ell -
1$) gives:
$$ |\calA| \leq 2 \cdot \left(\frac{3}{4}+o(1)\right)\frac{2^k}{\sqrt{\pi\ell}} =
\left(\frac{3}{2}+o(1)\right)\frac{2^k}{\sqrt{\pi\ell}}~,$$ and as
$|\calB_0| \leq 2^{h-m-1}$, \eqref{eq-B0-suff-condition} holds for
$\alpha=\frac{3}{4}$, as required.
\end{proof}

This completes the proof of Claim \ref{clm-rA=Omega(n)} and of Lemma
\ref{lem-rA+sA=Omega(n)}.

\section{Concluding remarks and open
problems}\label{sec::concluding}
\begin{itemize}
 \item We have shown that if two families of subsets of an $n$-element set, $\calA,\calB$, are
  $\ell$-cross-intersecting, and $\ell$ is sufficiently large, then $|\calA||\calB| \leq
  \binom{2\ell}{\ell}2^{n-2\ell}$, and in addition, we have given a complete characterization of all
  the extremal pairs $\calA,\calB$ for which equality is achieved.
 \item It would be interesting to prove that the above result holds
 for all values of $\ell$ (instead of all $\ell \geq \ell_0$ for
 some $\ell_0$). Perhaps knowing the precise structure of the
 extremal pairs $\calA,\calB$, as described in Theorem \ref{thm-1} (assuming that this holds for all $\ell$), will assist in proving this result.
 \item Finally, one may consider the corresponding problem where the pair $\calA,\calB$ does not have
 one possible cross-intersection, but rather a set $L$ of legal cross-intersections. Such notions have been
 studied in \cite{ACZ}, \cite{Sgall}, \cite{KeevashSudakov}, with different restrictions on $L$, and it
 would be interesting to derive tight bounds on $|\calA||\calB|$,
 and possibly describe the structure of all the extremal pairs, when
in addition, each member of $L$ is larger than some predefined
integer $\ell$.
\end{itemize}

\noindent\textbf{Acknowledgement} The authors wish to thank Benny
Sudakov for useful discussions.

\end{document}